\catcode`@=11
\expandafter\ifx\csname capitolo\endcsname\relax
      {\message {lettura formati e definizioni
standard}}\else{\endinput}\fi
\magnification \magstep1
\tolerance=1500\frenchspacing
%
% Fonti di caratteri
%

\font\titfnt=cmbx10 scaled \magstep1
\font\pagfnt=cmsl10
\font\parfnt=cmbx10 scaled \magstep1
\font\sprfnt=cmsl10 scaled \magstep1
\font\teofnt=cmssbx10
\font\tenenufnt=cmsl10
\font\nineenufnt=cmsl9
\font\eightenufnt=cmsl8

\font\ninerm=cmr9
\font\eightrm=cmr8
\font\sixrm=cmr6
 
\font\ninei=cmmi9
\font\eighti=cmmi8
\font\sixi=cmmi6
\skewchar\ninei='177 \skewchar\eighti='177
\skewchar\sixi='177
 
\font\ninesy=cmsy9
\font\eightsy=cmsy8
\font\sixsy=cmsy6
\skewchar\ninesy='60 \skewchar\eightsy='60
\skewchar\sixsy='60

\font\ninebf=cmbx9
\font\eightbf=cmbx8
\font\sixbf=cmbx6
\font\fivebf=cmbx5
 
\font\ninett=cmtt9
\font\eighttt=cmtt8
 
\hyphenchar\tentt=-1 % inhibit hyphenation in typewriter type
\hyphenchar\ninett=-1
\hyphenchar\eighttt=-1
 
\font\ninesl=cmsl9
\font\eightsl=cmsl8
 
\font\nineit=cmti9
\font\eightit=cmti8
 
\font\tenmat=cmss10

\font\sevenmat=cmss7 
 
\font\fivemat=cmss5
\newfam\matfam

\newskip\ttglue
%
% Definizione formato della pagina
%
\newdimen\pagewidth \newdimen\pageheight
\newdimen\ruleht
\hsize=31.6pc  \vsize=44.5pc  \maxdepth=2.2pt
\parindent=19pt
\parindent=19pt
\hfuzz=2pt
\pagewidth=\hsize \pageheight=\vsize \ruleht=.5pt
\abovedisplayskip=6pt plus 3pt minus 1pt
\belowdisplayskip=6pt plus 3pt minus 1pt
\abovedisplayshortskip=0pt plus 3pt
\belowdisplayshortskip=4pt plus 3pt
\def\blank{\vskip 12pt}
\def\blankii{\blank\blank}

\def\blankq{\vskip 3pt}
%
% Definizione del corpo dei caratteri e dei formati correnti
%
\def\tenpoint{\def\rm{\fam0\tenrm}%
  \def\enufnt{\tenenufnt}%
  \textfont0=\tenrm \scriptfont0=\sevenrm
\scriptscriptfont0=\fiverm
  \textfont1=\teni \scriptfont1=\seveni
\scriptscriptfont1=\fivei
  \textfont2=\tensy \scriptfont2=\sevensy
\scriptscriptfont2=\fivesy
  \textfont3=\tenex \scriptfont3=\tenex
\scriptscriptfont3=\tenex
  \def\it{\fam\itfam\tenit}%
  \textfont\itfam=\tenit
  \def\sl{\fam\slfam\tensl}%
  \textfont\slfam=\tensl
  \def\bf{\fam\bffam\tenbf}%
  \textfont\bffam=\tenbf \scriptfont\bffam=\sevenbf
   \scriptscriptfont\bffam=\fivebf
 \def\mat{\fam\matfam\tenmat} \textfont\matfam=\tenmat
     \scriptfont\matfam=\sevenmat \scriptscriptfont\matfam=\fivemat
  \def\tt{\fam\ttfam\tentt}%
  \textfont\ttfam=\tentt
  \tt \ttglue=.5em plus.25em minus.15em
  \normalbaselineskip=12pt
  \let\sc=\eightrm
  \let\big=\tenbig
  \setbox\strutbox=\hbox{\vrule height8.5pt depth3.5pt width\z@}%
  \normalbaselines\rm}
 
\def\ninepoint{\def\rm{\fam0\ninerm}%
  \def\enufnt{\nineenufnt}%
  \textfont0=\ninerm \scriptfont0=\sixrm
\scriptscriptfont0=\fiverm
  \textfont1=\ninei \scriptfont1=\sixi
\scriptscriptfont1=\fivei
  \textfont2=\ninesy \scriptfont2=\sixsy
\scriptscriptfont2=\fivesy
  \textfont3=\tenex \scriptfont3=\tenex
\scriptscriptfont3=\tenex
  \def\it{\fam\itfam\nineit}%
  \textfont\itfam=\nineit
  \def\sl{\fam\slfam\ninesl}%
  \textfont\slfam=\ninesl
  \def\bf{\fam\bffam\ninebf}%
  \textfont\bffam=\ninebf \scriptfont\bffam=\sixbf
   \scriptscriptfont\bffam=\fivebf
  \def\tt{\fam\ttfam\ninett}%
  \textfont\ttfam=\ninett
  \tt \ttglue=.5em plus.25em minus.15em
  \normalbaselineskip=11pt
  \let\sc=\sevenrm
  \let\big=\ninebig
  \setbox\strutbox=\hbox{\vrule height8pt depth3pt width\z@}%
  \normalbaselines\rm}
 
\def\eightpoint{\def\rm{\fam0\eightrm}%
  \def\enufnt{\eightenufnt}%
  \textfont0=\eightrm \scriptfont0=\sixrm
\scriptscriptfont0=\fiverm
  \textfont1=\eighti \scriptfont1=\sixi
\scriptscriptfont1=\fivei
  \textfont2=\eightsy \scriptfont2=\sixsy
\scriptscriptfont2=\fivesy
  \textfont3=\tenex \scriptfont3=\tenex
\scriptscriptfont3=\tenex
  \def\it{\fam\itfam\eightit}%
  \textfont\itfam=\eightit
  \def\sl{\fam\slfam\eightsl}%
  \textfont\slfam=\eightsl
  \def\bf{\fam\bffam\eightbf}%
  \textfont\bffam=\eightbf \scriptfont\bffam=\sixbf
   \scriptscriptfont\bffam=\fivebf
  \def\tt{\fam\ttfam\eighttt}%
  \textfont\ttfam=\eighttt
  \tt \ttglue=.5em plus.25em minus.15em
  \normalbaselineskip=9pt
  \let\sc=\sixrm
  \let\big=\eightbig
  \setbox\strutbox=\hbox{\vrule height7pt depth2pt width\z@}%
  \normalbaselines\rm}
 
 % use after $ in ninepoint sections
\def\tenbig#1{{\hbox{$\left#1\vbox to8.5pt{}\right.\n@space$}}}
\def\ninebig#1{{\hbox{$\textfont0=\tenrm\textfont2=\tensy
  \left#1\vbox to7.25pt{}\right.\n@space$}}}
\def\eightbig#1{{\hbox{$\textfont0=\ninerm\textfont2=\ninesy
  \left#1\vbox to6.5pt{}\right.\n@space$}}}
%
% Definizione macro di allineamento testi
%
\def\corsivo{\enufnt}

\def\footnote#1{\edef\@sf{\spacefactor\the\spacefactor}$^{#1}$\@sf
      \insert\footins\bgroup\ninepoint
      \interlinepenalty100 \let\par=\endgraf
        \leftskip=\z@skip \rightskip=\z@skip
        \splittopskip=10pt plus 1pt minus 1pt \floatingpenalty=20000
        \smallskip\item{$^{#1}$}\bgroup\strut\aftergroup\@foot\let\next}
\skip\footins=12pt plus 2pt minus 4pt % space added when footnote is present
%\count\footins=1000 % footnote magnification factor (1 to 1)
\dimen\footins=30pc % maximum footnotes per page
%
%  Registrazione referenze su file esterno e macro di conteggio
%
\newcount\numbibliogr@fi@
\global\numbibliogr@fi@=1
\newwrite\fileref
\immediate\openout\fileref=\jobname.ref
\immediate\write\fileref{\parindent 30pt}
\def\cita#1#2{\def\us@gett@{\the\numbibliogr@fi@}
\expandafter\xdef\csname:bib_#1\endcsname{\us@gett@}
           \immediate\write\fileausiliario{\string\expandafter\string\edef
           \string\csname:bib_#1\string\endcsname{\us@gett@}}
           \immediate\write\fileref
           {\par\noexpand\item{{[\the\numbibliogr@fi@]\enspace}}}\ignorespaces
           \immediate\write\fileref{{#2}}\ignorespaces
           \global\advance\numbibliogr@fi@ by 1\ignorespaces}
\def\bibref#1{\seindefinito{:bib_#1}
          \immediate\write16{ !!! \string\bibref{#1} non definita !!!}
 
\expandafter\xdef\csname:bib_#1\endcsname{??}\fi
      {$^{[\csname:bib_#1\endcsname]}$}}
\def\dbiref#1{\seindefinito{:bib_#1}
          \immediate\write16{ !!! \string\bibref{#1} non definita !!!}
 
\expandafter\xdef\csname:bib_#1\endcsname{??}\fi
      {[\csname:bib_#1\endcsname]}}

\def\references{\immediate\closeout\fileref
                \par\goodbreak
                \blankii
                \centerline{\parfnt References}
                \nobreak\blank\nobreak
                \input \jobname.ref}

%
% Definizione macro di intestazione
%
\def\title#1{\null\blankii\noindent{\titfnt\uppercase{#1}}\blank}

\def\author#1{\leftskip 1.8cm\smallskip\noindent{#1}\smallskip\leftskip 0pt}
\def\abstract#1{\par\blankii\noindent
          {\ninepoint
            {\bf Abstract. }{\rm #1}}
          \par}

\def\sunto#1{\par\blankii\noindent
          {\ninepoint
            {\bf Sunto. }{\rm #1}}
          \par}
%
% Definizione macro per la numerazione delle formule, capitoli,
% paragrafi, ecc...
%
%
%  Definizioni generali
\def\seindefinito#1{\expandafter\ifx\csname#1\endcsname\relax}
\newwrite\filesimboli
    \immediate\openout\filesimboli=\jobname.smb
\newwrite\fileausiliario
    \openin1 \jobname.aux
 
\ifeof1\relax\else\closein1\relax\input\jobname.aux\fi
    \immediate\openout\fileausiliario=\jobname.aux
%
%  Paragrafo
\newdimen\@mpiezz@
\@mpiezz@=\hsize
\newbox\boxp@r@gr@fo
\newcount\nump@r@gr@fo
\global\nump@r@gr@fo=0
\def\section#1#2{
           \global\advance\nump@r@gr@fo by 1
           \xdef\paragrafocorrente{\the\nump@r@gr@fo}
\expandafter\xdef\csname:sec_#1\endcsname{\paragrafocorrente}
           \global\numsottop@r@gr@fo=0
      \immediate\write\filesimboli{  }
      \immediate\write\filesimboli{  }
      \immediate\write\filesimboli{--> Section \paragrafocorrente :
               rif.: #1}
\immediate\write\fileausiliario{\string\expandafter\string\edef
\string\csname:sec_#1\string\endcsname{\paragrafocorrente}}
      \write16{Section {\csname:sec_#1\endcsname}.}
           \goodbreak\vskip 24pt plus 6pt\noindent\ignorespaces
              {\setbox\boxp@r@gr@fo=\hbox{\parfnt\noindent\ignorespaces
              {\paragrafocorrente.\quad}}\ignorespaces
             \advance\@mpiezz@ by -\wd\boxp@r@gr@fo
             \box\boxp@r@gr@fo\vtop{\hsize=\@mpiezz@\noindent\parfnt #2}}
           \par\nobreak\vskip 12pt plus 3pt\nobreak
           \noindent\ignorespaces}
\def\secref#1{\seindefinito{:sec_#1}
          \immediate\write16{ !!! \string\secref{#1} non definita !!!}
 \expandafter\xdef\csname:sec_#1\endcsname{??}\fi
      \csname:sec_#1\endcsname
      }
%  Appendice
\newcount\num@ppendice
\global\num@ppendice=64
\def\appendix#1#2{
           \global\advance\num@ppendice by 1
           \xdef\paragrafocorrente{\char\the\num@ppendice}
           \expandafter\xdef\csname:app_#1\endcsname{\paragrafocorrente}
           \global\numsottop@r@gr@fo=0
      \immediate\write\filesimboli{  }
      \immediate\write\filesimboli{  }
      \immediate\write\filesimboli{--> Appendice \paragrafocorrente :
               rif.: #1}
      \immediate\write\fileausiliario{\string\expandafter\string\edef
      \string\csname:app_#1\string\endcsname{\paragrafocorrente}}
      \write16{Appendice {\csname:app_#1\endcsname}.}
           \goodbreak\vskip 18pt plus 6pt\noindent\ignorespaces
              {\setbox\boxp@r@gr@fo=\hbox{\parfnt\noindent\ignorespaces
              {\paragrafocorrente.\quad}}\ignorespaces
             \advance\@mpiezz@ by -\wd\boxp@r@gr@fo
             \box\boxp@r@gr@fo\vtop{\hsize=\@mpiezz@\noindent\parfnt #2}}
           \par\nobreak\vskip 9pt plus 3pt\nobreak
           \noindent\ignorespaces}
\def\appref#1{\seindefinito{:app_#1}
          \immediate\write16{ !!! \string\appref{#1} non definita !!!}
          \expandafter\xdef\csname:app_#1\endcsname{??}\fi
      \csname:app_#1\endcsname
      }
%
%  Sottoparagrafo
\newcount\numsottop@r@gr@fo
\global\numsottop@r@gr@fo=0
\def\subsection#1#2{
           \global\advance\numsottop@r@gr@fo by 1
 \xdef\sottoparcorrente{\paragrafocorrente.\the\numsottop@r@gr@fo}
 \expandafter\xdef\csname:sbs_#1\endcsname{\sottoparcorrente}
      \immediate\write\filesimboli{  }
      \immediate\write\filesimboli{-----> Subsection \sottoparcorrente :
               rif.: #1}
 \immediate\write\fileausiliario{\string\expandafter\string\edef
 \string\csname:sbs_#1\string\endcsname{\sottoparcorrente}}
      \write16{Subsection {\csname:sbs_#1\endcsname}.}
           \goodbreak\vskip 9pt plus 2pt\noindent\ignorespaces
 {\setbox\boxp@r@gr@fo=\hbox{\sprfnt\noindent\ignorespaces
              {\sottoparcorrente\quad}}\ignorespaces
             \advance\@mpiezz@ by -\wd\boxp@r@gr@fo
 \box\boxp@r@gr@fo\vtop{\hsize=\@mpiezz@\noindent\sprfnt#2}}
           \par\nobreak\vskip 3pt plus 1pt\nobreak
           \noindent\ignorespaces}
\def\sbsref#1{\seindefinito{:sbs_#1}
          \immediate\write16{ !!! \string\sbsref{#1} non definita !!!}
 \expandafter\xdef\csname:sbs_#1\endcsname{??}\fi
      \csname:sbs_#1\endcsname
      }
%
%  Formula
\def\eqalignno#1{\leqalignno{#1}}
\newcount\numformul@
\global\numformul@=0
\def\formula#1{
            \formdef{#1}\ignorespaces
            \leqno{(\csname:frm_#1\endcsname)}
            }
\def\formdef#1{
            \global\advance\numformul@ by 1
            \def\us@gett@{\the\numformul@}
 \expandafter\xdef\csname:frm_#1\endcsname{\us@gett@}
      \immediate\write\filesimboli{--------> Formula \us@gett@ :
               rif.: #1}
 \immediate\write\fileausiliario{\string\expandafter\string\edef
 \string\csname:frm_#1\string\endcsname{\us@gett@}}
            }
\def\frmref#1{\seindefinito{:frm_#1}
          \immediate\write16{ !!! \string\frmref{#1} non definita !!!}
 \expandafter\xdef\csname:frm_#1\endcsname{??}\fi
      (\csname:frm_#1\endcsname)}
%
%  Teorema
\newcount\numenunci@to
\global\numenunci@to=0
\def\endclaim{\endgroup
            \par\if F\sp@zi@tur@{\blankq}\gdef\sp@zi@tur@{T}\fi}
\newcount\numteorema
\global\numteorema=0
\def\theorem#1{
            \global\advance\numteorema by 1
            \def\us@gett@{\the\numteorema}
 \expandafter\xdef\csname:thr_#1\endcsname{\us@gett@}
      \immediate\write\filesimboli{--------> Theorem \us@gett@ :
               rif.: #1}
 \immediate\write\fileausiliario{\string\expandafter\string\edef
 \string\csname:thr_#1\string\endcsname{\us@gett@}}
            \par\if T\sp@zi@tur@{\gdef\sp@zi@tur@{F}}\else{\blankq}\fi
            \noindent{\teofnt Theorem \us@gett@:\quad}\begingroup\enufnt
            \ignorespaces}
\def\theoremnn{
            \par\if T\sp@zi@tur@{\gdef\sp@zi@tur@{F}}\else{\blankq}\fi
            \noindent{\teofnt Theorem:\quad}\begingroup\enufnt
            \ignorespaces}
\def\theoremtx#1#2{
            \global\advance\numteorema by 1
            \def\us@gett@{\the\numteorema}
 \expandafter\xdef\csname:thr_#1\endcsname{\us@gett@}
      \immediate\write\filesimboli{--------> Theorem \us@gett@ :
               rif.: #1}
 \immediate\write\fileausiliario{\string\expandafter\string\edef
 \string\csname:thr_#1\string\endcsname{\us@gett@}}
            \par\if T\sp@zi@tur@{\gdef\sp@zi@tur@{F}}\else{\blankq}\fi
            \noindent{\teofnt Theorem \us@gett@:
                {\enufnt #2}.\quad}\begingroup\enufnt
            \ignorespaces}
\def\thrref#1{\seindefinito{:thr_#1}
          \immediate\write16{ !!! \string\thrref{#1} non definita !!!}
 \expandafter\xdef\csname:thr_#1\endcsname{??}\fi
      \csname:thr_#1\endcsname}
%
%  Proposizione
\newcount\numproposizione
\global\numproposizione=0
\def\proposition#1{
            \global\advance\numproposizione by 1
            \def\us@gett@{\the\numproposizione}
 \expandafter\xdef\csname:pro_#1\endcsname{\us@gett@}
      \immediate\write\filesimboli{--------> Proposition \us@gett@ :
               rif.: #1}
 \immediate\write\fileausiliario{\string\expandafter\string\edef
 \string\csname:pro_#1\string\endcsname{\us@gett@}}
            \par\if T\sp@zi@tur@{\gdef\sp@zi@tur@{F}}\else{\blankq}\fi
            \noindent{\teofnt Proposition \us@gett@:\quad}\begingroup\enufnt
            \ignorespaces}
\def\propositionn{
            \par\if T\sp@zi@tur@{\gdef\sp@zi@tur@{F}}\else{\blankq}\fi
            \noindent{\teofnt Proposition:\quad}\begingroup\enufnt
            \ignorespaces}
\def\propositiontx#1#2{
            \global\advance\numproposizione by 1
            \def\us@gett@{\the\numproposizione}
 \expandafter\xdef\csname:pro_#1\endcsname{\us@gett@}
      \immediate\write\filesimboli{--------> Proposition \us@gett@ :
               rif.: #1}
  \immediate\write\fileausiliario{\string\expandafter\string\edef
 \string\csname:pro_#1\string\endcsname{\us@gett@}}
            \par\if T\sp@zi@tur@{\gdef\sp@zi@tur@{F}}\else{\blankq}\fi
            \noindent{\teofnt Proposition \us@gett@:
               {\enufnt #2}.\quad}\begingroup\enufnt
            \ignorespaces}
\def\proref#1{\seindefinito{:pro_#1}
          \immediate\write16{ !!! \string\proref{#1} non definita !!!}
 \expandafter\xdef\csname:pro_#1\endcsname{??}\fi
      \csname:pro_#1\endcsname}
%
%  Corollario
\newcount\numcorollario
\global\numcorollario=0
\def\corollary#1{
            \global\advance\numcorollario by 1
            \def\us@gett@{\the\numcorollario}
 \expandafter\xdef\csname:cor_#1\endcsname{\us@gett@}
      \immediate\write\filesimboli{--------> Corollary \us@gett@ :
               rif.: #1}
 \immediate\write\fileausiliario{\string\expandafter\string\edef
 \string\csname:cor_#1\string\endcsname{\us@gett@}}
            \par\if T\sp@zi@tur@{\gdef\sp@zi@tur@{F}}\else{\blankq}\fi
            \noindent{\teofnt Corollary
\us@gett@:\quad}\begingroup\enufnt
            \ignorespaces}
\def\corref#1{\seindefinito{:cor_#1}
          \immediate\write16{ !!! \string\corref{#1} non definita !!!}
 \expandafter\xdef\csname:cor_#1\endcsname{??}\fi
      \csname:cor_#1\endcsname}
%
%  Lemma
\newcount\numlemma
\global\numlemma=0
\def\lemma#1{
            \global\advance\numlemma by 1
            \def\us@gett@{\the\numlemma}
\expandafter\xdef\csname:lem_#1\endcsname{\us@gett@}
      \immediate\write\filesimboli{--------> Lemma \us@gett@ :
               rif.: #1}
 \immediate\write\fileausiliario{\string\expandafter\string\edef
 \string\csname:lem_#1\string\endcsname{\us@gett@}}
            \par\if T\sp@zi@tur@{\gdef\sp@zi@tur@{F}}\else{\blankq}\fi
            \noindent{\teofnt Lemma \us@gett@:\quad}\begingroup\enufnt
            \ignorespaces}
\def\lemref#1{\seindefinito{:lem_#1}
          \immediate\write16{ !!! \string\lemref{#1} non definita !!!}
 \expandafter\xdef\csname:lem_#1\endcsname{??}\fi
      \csname:lem_#1\endcsname}
%
%  Definizione
\newcount\numdefinizione
\global\numdefinizione=0
\def\definition#1{
            \global\advance\numdefinizione by 1
            \def\us@gett@{\the\numdefinizione}
 \expandafter\xdef\csname:def_#1\endcsname{\us@gett@}
      \immediate\write\filesimboli{--------> Definition \us@gett@ :
               rif.: #1}
 \immediate\write\fileausiliario{\string\expandafter\string\edef
 \string\csname:def_#1\string\endcsname{\us@gett@}}
            \par\if T\sp@zi@tur@{\gdef\sp@zi@tur@{F}}\else{\blankq}\fi
            \noindent{\teofnt Definition \us@gett@:\quad}\begingroup\enufnt
            \ignorespaces}
\def\definitiontx#1#2{
            \global\advance\numdefinizione by 1
            \def\us@gett@{\the\numdefinizione}
 \expandafter\xdef\csname:def_#1\endcsname{\us@gett@}
      \immediate\write\filesimboli{--------> Definition \us@gett@ :
               rif.: #1}
 \immediate\write\fileausiliario{\string\expandafter\string\edef
 \string\csname:def_#1\string\endcsname{\us@gett@}}
            \par\if T\sp@zi@tur@{\gdef\sp@zi@tur@{F}}\else{\blankq}\fi
            \noindent{\teofnt Definition \us@gett@:
               {\enufnt #2}.\quad}\begingroup\enufnt
            \ignorespaces}
\def\defref#1{\seindefinito{:def_#1}
          \immediate\write16{ !!! \string\defref{#1} non definita !!!}
 \expandafter\xdef\csname:def_#1\endcsname{??}\fi
      \csname:def_#1\endcsname}
%
%  Dimostrazione
\def\proof{\par\if T\sp@zi@tur@{\gdef\sp@zi@tur@{F}}\else{\blankq}\fi
    \noindent{\teofnt Proof.\quad}\begingroup\ignorespaces}
\def\prooftx#1{\par\if T\sp@zi@tur@{\gdef\sp@zi@tur@{F}}\else{\blankq}\fi
    \noindent{\teofnt Proof #1.\quad}\begingroup\ignorespaces}
\def\endproof{\nobreak\quad\nobreak\hfill\nobreak{\enufnt Q.E.D.}
    \endclaim}
%
% Figura
\newcount\numfigur@
\global\numfigur@=0
\newbox\boxfigur@
\newbox\comfigur@
\newdimen\@mpfigur@
\newdimen\m@rfigur@
%\m@rfigur@=.4 cm
\m@rfigur@=2 pc
\@mpfigur@=\hsize
\advance\@mpfigur@ by -2\m@rfigur@
\def\figure#1#2#3{
            \global\advance\numfigur@ by 1
            \def\us@gett@{\the\numfigur@}
            \expandafter\xdef\csname:fig_#1\endcsname{\us@gett@}
      \immediate\write\filesimboli{--------> Figure \us@gett@ :
               rif.: #1}
      \immediate\write\fileausiliario{\string\expandafter\string\edef
      \string\csname:fig_#1\string\endcsname{\us@gett@}}
      \setbox\boxfigur@\vbox{\centerline{#2}}
      \topinsert
         {\vbox{
               \vskip 1pt
               \box\boxfigur@
               \vskip 4pt}}
      \setbox\comfigur@\vtop{\hsize=\@mpfigur@\parindent 0pt
         {\ninepoint
         {\teofnt Figure \us@gett@.}\enspace{#3}}
      }
      \centerline{\box\comfigur@}
      \endinsert
      \write16{Figure {\csname:fig_#1\endcsname}.}
      }
\def\figcont#1#2{
      \setbox\boxfigur@\vbox{\centerline{#2}}
      \topinsert
         {\vbox{
               \vskip 1pt
               \box\boxfigur@
               \vskip 4pt}}
 \setbox\comfigur@\vtop{\hsize=\@mpfigur@\parindent 0pt
         {\ninepoint
         {\teofnt Figure \figref{#1}.}\enspace{(continued).}}
      }
      \centerline{\box\comfigur@}
      \endinsert
      \write16{Figura (cont) {\csname:fig_#1\endcsname}.}
}
\def\figref#1{\seindefinito{:fig_#1}
          \immediate\write16{ !!! \string\figref{#1} non definita !!!}
 \expandafter\xdef\csname:fig_#1\endcsname{??}\fi
      \csname:fig_#1\endcsname}
% Citazione
\def\citazione{\begingroup
     \everypar{\parshape 1 \m@rfigur@ \@mpfigur@\noindent}}

%
%  Osservazione/i
\def\remarks{\par\if
T\sp@zi@tur@{\gdef\sp@zi@tur@{F}}\else{\blankq}\fi
    \noindent{\teofnt Remarks.\quad}\ignorespaces}
\def\remark{\par\if
T\sp@zi@tur@{\gdef\sp@zi@tur@{F}}\else{\blankq}\fi
    \noindent{\teofnt Remark.\quad}\ignorespaces}
%
% Tavola
% #1 : numerazione
% #2 : didascalia
% #3 : tabella
\newcount\numt@vol@
\global\numt@vol@=0
\newbox\boxt@vol@
\newbox\comt@vol@
\newdimen\@mpt@vol@
\newdimen\m@rt@vol@
%\m@rt@vol@=.4 cm
\m@rt@vol@=2 pc
\@mpt@vol@=\hsize
\advance\@mpt@vol@ by -2\m@rt@vol@
\def\table#1#2#3{
            \global\advance\numt@vol@ by 1
            \def\us@gett@{\the\numt@vol@}
            \expandafter\xdef\csname:tav_#1\endcsname{\us@gett@}
      \immediate\write\filesimboli{--------> Table \us@gett@ :
               rif.: #1}
      \immediate\write\fileausiliario{\string\expandafter\string\edef
      \string\csname:tav_#1\string\endcsname{\us@gett@}}
      \setbox\boxt@vol@\vbox{#3}
      \topinsert
      \setbox\comt@vol@\vtop{\hsize=\@mpt@vol@\parindent 0pt
         {\ninepoint
         {\teofnt Table \us@gett@.}\enspace{#2}}
      }
      \centerline{\box\comt@vol@}
      \vskip 5pt%
      \centerline{\box\boxt@vol@}%
      \endinsert
      \write16{Table {\csname:tav_#1\endcsname}.}
}
\def\tabcont#1#2{
      \setbox\boxt@vol@\vbox{#2}
      \topinsert
      \setbox\comt@vol@\vtop{\hsize=\@mpt@vol@\parindent 0pt
         {\ninepoint
         {\teofnt Table \tabref{#1}.}\enspace{(continued)}}
      }
      \centerline{\box\comt@vol@}
      \vskip 5pt%
      \centerline{\box\boxt@vol@}%
      \endinsert
      \write16{Table {\csname:tav_#1\endcsname}.}
}
\def\tabref#1{\seindefinito{:tav_#1}
          \immediate\write16{ !!! \string\tabref{#1} non definita !!!}
          \expandafter\xdef\csname:tav_#1\endcsname{??}\fi
      \csname:tav_#1\endcsname}
%
% Inibizione spaziatura prima degli enunciati
%

\def\noblank{\gdef\sp@zi@tur@{T}}
\def\incolonna#1{\displ@y \tabskip=\centering
   \halign to \displaywidth{\hfil$\@lign\displaystyle{{}##}$\tabskip=0pt
       &$\@lign\displaystyle{{}##}$\hfil\tabskip=\centering
       &\llap{$\@lign\displaystyle{{}##}$}\tabskip=0pt\crcr
       #1\crcr}}
%
% Definizioni per l'output
%
\nopagenumbers
\def\testos{\null}
\def\testod{\null}
\headline={\if T\tpage{\gdef\tpage{F}{\hfil}}
 \else{\ifodd\pageno\rightheadline\else\leftheadline\fi}
           \fi}
 
\gdef\tpage{T}
\def\rightheadline{\hfil{\pagfnt\testod}\hfil{\pagfnt\folio}}
\def\leftheadline{{\pagfnt\folio}\hfil{\pagfnt\testos}\hfil}
\voffset=2\baselineskip
\everypar={\gdef\sp@zi@tur@{F}}
\catcode`@=12
%
%   Definizioni simboli d'uso comune

\def\pmb#1{\setbox0=\hbox{#1}\ignorespaces
    \hbox{\kern-.02em\copy0\kern-\wd0\ignorespaces
    \kern.05em\copy0\kern-\wd0\ignorespaces
    \kern-.02em\raise.02em\box0 }}
\def\gt{>}
\def\lt{<}
\def\rho{\varrho}
\def\theta{\vartheta}
       % questa istruzione DEVE precedere la
                       % ridefinizione di \phi della prossima riga
\def\phi{\varphi}
\def\epsilon{\varepsilon}

\def\frac#1#2{{{#1}\over{#2}}}

\def\reali{\mathinner{\bf R}}
\def\complessi{\mathinner{\bf C}}

\def\interi{\mathinner{\bf Z}}

\def\poisson#1#2{\lbrace#1,#2\rbrace}

\tenpoint\rm
%DEFINIZIONI LOCALI

\def\epsilon{\varepsilon}
\def\phi{\varphi}

\def\rho{\varrho}
\def\theta{\vartheta}

\def\lie#1{L_{#1}}

\def\Cscr{{\cal C}}

\def\Jscr{{\cal J}}
\def\Kscr{{\cal K}}

\def\Nscr{{\cal N}}
\def\Pscr{{\cal P}}

\def\Rscr{{\cal R}}
\def\Wscr{{\cal W}}
\def\Zscr{{\cal Z}}
\def\oxi{{\mathaccent"17 \xi}}
\def\oeta{{\mathaccent"17 \eta}}
\def\ozeta{{\mathaccent"17 \zeta}}

\def\poisson#1#2{\lbrace #1,#2 \rbrace}

\def\lgd{\mathop{{\rm log}_2}}
\def\span{\mathop {\rm span}}
\def\eqalignno#1{\leqalignno{#1}}
\def\norma#1{\left\|{#1}\right\|}
\def\ordnorma#1{\|{#1}\|}
\def\bignorma#1{\bigl\|{#1}\bigr\|}

%FINE DEFINIZIONI LOCALI
%BIBLIOGRAFIA

\cita{Birkhoff-27}{Birkhoff, G. D.: {\it Dynamical systems}, New York 
(1927).}

\cita{Cherry-1926}{Cherry, T. M.: {\it On the solutions of Hamiltonian
systems in the neighborhood of a singular point}, Proc. London
Math. Soc., Ser. 2, {\bf 27}, 151--170 (1926).}

\cita{Giorgilli-95}{Giorgilli, A.: {\it Quantitative methods in classical
perturbation theory}, proceedings of the Nato ASI school ``From Newton to
chaos: modern techniques for understanding and coping with chaos in
N--body dynamical systems'', A.E.{\ }Roy e B.D.{\ }Steves eds., Plenum Press, 
New York (1995).} 

\cita{Giorgilli-97.3}{Giorgilli, A.  and Locatelli, U.: {\it Kolmogorov theorem 
and classical perturbation theory}, ZAMP {\bf 48}, 220--261 (1997).}

\cita{Giorgilli-97.4}{Giorgilli, A.  and Locatelli, U.: {\it On classical series
expansions for quasi--periodic motions}, MPEJ {\bf 3} N. 5 (1997).}

\cita{Giorgilli-98a}{Giorgilli: {\it Classical constructive 
methods in KAM theory}, PSS, {\bf 46} 1441--1451 (1998).}

\cita{Giorgilli-99}{Giorgilli, A.  and Locatelli, U.: {\it A classical
self--contained proof of Kolmogorov's theorem on invariant tori}, in
{\it Hamiltonian systems with three or more degrees of freedom},
Carles Sim\'o ed., NATO ASI series C, Vol. 533, Kluwer Academic
Publishers, Dordrecht--Boston--London (1999).}

\cita{Giorgilli-2001}{Giorgilli, A.: {\it Unstable equilibria of
Hamiltonian systems}, Disc. and Cont. Dynamical Systems, Vol. 7, N. 4,
855--871 (2001).}

\cita{Groebner-67}{Gr\"obner, W.: {\it Die  Lie--Reihen und Ihre 
Anwendungen}, VEB Deutscher Verlag der Wissenschaften (1967).}

\cita{Lyapounov-1892}{Lyapunov, A.M.: {\it The General Problem of the
Stability of Motion} (In Russian), Doctoral dissertation,
Univ. Kharkov (1892).  French translation in: {\it Probl\`eme
g\'en\'eral de la stabilit\'e du mouvement}, Annales de la Facult\'e
des Sciences de Toulouse, deuxi\`eme s\'erie, Tome IX, 203--474
(1907). Reprinted in: Ann. Math. Study, Princeton University Press,
n.\ 17, (1949).}

\cita{Moser-1958}{Moser, J.: {\it On the generalization of a theorem of
A. Liapounoff}, Comm. Pure Appl. Math. {\bf 11}, 257--271 (1958).}

\cita{Siegel-1971}{Siegel, C.L. and Moser, J.K.: {\it Lectures in
Celestial Mechanics,} Springer--Verlag, Berlin Heidelberg New York (1971).}

%FINE BIBLIOGRAFIA
\expandafter\edef\csname:bib_Birkhoff-27\endcsname{1}
\expandafter\edef\csname:bib_Cherry-1926\endcsname{2}
\expandafter\edef\csname:bib_Giorgilli-95\endcsname{3}
\expandafter\edef\csname:bib_Giorgilli-97.3\endcsname{4}
\expandafter\edef\csname:bib_Giorgilli-97.4\endcsname{5}
\expandafter\edef\csname:bib_Giorgilli-98a\endcsname{6}
\expandafter\edef\csname:bib_Giorgilli-99\endcsname{7}
\expandafter\edef\csname:bib_Giorgilli-2001\endcsname{8}
\expandafter\edef\csname:bib_Groebner-67\endcsname{9}
\expandafter\edef\csname:bib_Lyapounov-1892\endcsname{10}
\expandafter\edef\csname:bib_Moser-1958\endcsname{11}
\expandafter\edef\csname:bib_Siegel-1971\endcsname{12}
\expandafter\edef\csname:sec_1\endcsname{1}
\expandafter\edef\csname:frm_ham.1\endcsname{1}
\expandafter\edef\csname:frm_ham.2\endcsname{2}
\expandafter\edef\csname:frm_meln.1\endcsname{3}
\expandafter\edef\csname:frm_meln.2\endcsname{4}
\expandafter\edef\csname:frm_1.1\endcsname{5}
\expandafter\edef\csname:frm_1.2\endcsname{6}
\expandafter\edef\csname:thr_nrm.1\endcsname{1}
\expandafter\edef\csname:frm_lyap.1\endcsname{7}
\expandafter\edef\csname:thr_nrm.1.a\endcsname{2}
\expandafter\edef\csname:frm_lyap.1\endcsname{8}
\expandafter\edef\csname:sec_2\endcsname{2}
\expandafter\edef\csname:sbs_2.1\endcsname{2.1}
\expandafter\edef\csname:frm_alg.1\endcsname{9}
\expandafter\edef\csname:frm_alg.2\endcsname{10}
\expandafter\edef\csname:frm_alg.3\endcsname{11}
\expandafter\edef\csname:frm_nucleo-range\endcsname{12}
\expandafter\edef\csname:frm_gen_soluzione\endcsname{13}
\expandafter\edef\csname:sbs_2.2\endcsname{2.2}
\expandafter\edef\csname:frm_indici\endcsname{14}
\expandafter\edef\csname:frm_sottospazi\endcsname{15}
\expandafter\edef\csname:frm_splitting\endcsname{16}
\expandafter\edef\csname:frm_nonrisonanza\endcsname{17}
\expandafter\edef\csname:frm_gen.1\endcsname{18}
\expandafter\edef\csname:sbs_2.3\endcsname{2.3}
\expandafter\edef\csname:frm_indici.a\endcsname{19}
\expandafter\edef\csname:frm_splitting.a\endcsname{20}
\expandafter\edef\csname:sec_3\endcsname{3}
\expandafter\edef\csname:frm_disk\endcsname{21}
\expandafter\edef\csname:frm_lambda.disk\endcsname{22}
\expandafter\edef\csname:frm_norm.pol\endcsname{23}
\expandafter\edef\csname:lem_generatrici\endcsname{1}
\expandafter\edef\csname:sbs_3.2\endcsname{3.1}
\expandafter\edef\csname:lem_art.1\endcsname{2}
\expandafter\edef\csname:cor_art.1a\endcsname{1}
\expandafter\edef\csname:sbs_3.1\endcsname{3.2}
\expandafter\edef\csname:frm_tabella_poisson\endcsname{24}
\expandafter\edef\csname:frm_stimagen\endcsname{25}
\expandafter\edef\csname:frm_cauchy.f\endcsname{26}
\expandafter\edef\csname:frm_cauchy.f-beq\endcsname{27}
\expandafter\edef\csname:frm_cauchy.Z-dies\endcsname{28}
\expandafter\edef\csname:sbs_3.3\endcsname{3.3}
\expandafter\edef\csname:frm_jset\endcsname{29}
\expandafter\edef\csname:frm_T\endcsname{30}
\expandafter\edef\csname:frm_T-prop.1\endcsname{31}
\expandafter\edef\csname:frm_T-prop.2\endcsname{32}
\expandafter\edef\csname:frm_museq\endcsname{33}
\expandafter\edef\csname:lem_nrm.1\endcsname{3}
\expandafter\edef\csname:frm_nrm.1\endcsname{34}
\expandafter\edef\csname:frm_nrm.2\endcsname{35}
\expandafter\edef\csname:frm_nrm.2a\endcsname{36}
\expandafter\edef\csname:frm_nrm.3\endcsname{37}
\expandafter\edef\csname:frm_nrm.3a\endcsname{38}
\expandafter\edef\csname:frm_nrm.4\endcsname{39}
\expandafter\edef\csname:frm_stima.Z\endcsname{40}
\expandafter\edef\csname:frm_stima.Z.sharp\endcsname{41}
\expandafter\edef\csname:frm_D\endcsname{42}
\expandafter\edef\csname:frm_sta.1\endcsname{43}
\expandafter\edef\csname:frm_lchi.Z.1\endcsname{44}
\expandafter\edef\csname:frm_lchi.Z.2\endcsname{45}
\expandafter\edef\csname:frm_lchi.H.1\endcsname{46}
\expandafter\edef\csname:frm_lchi.H.2\endcsname{47}
\expandafter\edef\csname:frm_museq.1\endcsname{48}
\expandafter\edef\csname:sbs_3.4\endcsname{3.4}
\expandafter\edef\csname:frm_nrm.11\endcsname{49}
\expandafter\edef\csname:frm_nuseq\endcsname{50}
\expandafter\edef\csname:frm_nrm.16\endcsname{51}
\expandafter\edef\csname:sec_4\endcsname{4}
\expandafter\edef\csname:app_lyap.app.a\endcsname{\char 65}
\expandafter\edef\csname:app_lyap.app.b\endcsname{\char 66}
\expandafter\edef\csname:frm_par_poisson\endcsname{52}
\expandafter\edef\csname:frm_ps.1\endcsname{53}
\expandafter\edef\csname:frm_ineq\endcsname{54}
\expandafter\edef\csname:frm_ps.3\endcsname{55}
\expandafter\edef\csname:frm_ps.2\endcsname{56}
\expandafter\edef\csname:frm_ps.5\endcsname{57}

\def\testos{A. Giorgilli}
\def\testod{On a theorem of Lyapounov}
\noindent
To appear in: \hfill\break
{\sl Rendiconti dell'Istituto Lombardo Accademia di Scienze e Lettere,
Classe di Scienze.}
\vskip 24 pt

\title{On a theorem of Lyapounov}

\author{\it ANTONIO GIORGILLI
\hfill\break Dipartimento di Matematica,
Via Saldini 50,
\hfill\break 20133\ ---\  Milano, Italy.}

\sunto{Si mostra che un sistema Hamiltoniano nell'intorno di un punto
di equilibrio, sotto condizione che gli autovalori soddisfino delle
condizioni di non--risonanza del tipo di Melnikov, ammette una forma
normale che rende evidente l'esistenza di una variet\`a invariante
(locale) a due dimensioni sulla quale si hanno soluzioni note.  Nel
caso di un autovalore puramente immaginario tali soluzioni formano una
famiglia periodica a due parametri che costituisce la continuazione
naturale di un modo normale.  Questo secondo risultato \`e stato
dimostrato in precedenza da Lyapounov.  In questo lavoro si completa
quello di Lyapounov dimostrando la convergenza della trasformazione
dell'Hamiltoniana a forma normale e rimuovendo le restrizione che gli
autovalori siano puramente immaginari.}

\abstract{It is shown that  a Hamiltonian system in the neighbourhood
of an equilibrium may be given a special normal form in case the
eigenvalues of the linearized system satisfy non--resonance conditions
of Melnikov's type.  The normal form possesses a two dimensional
(local) invariant manifold on which the solutions are known.  If the
eigenvalue is pure imaginary then these solutions are the natural
continuation of a normal mode of the linear system.  The latter result
was first proved by Lyapounov.  The present paper completes
Lyapounov's result in that the convergence of the transformation of
the Hamiltonian to a normal form is proven and the condition that the
eigenvalues be pure imaginary is removed.}

\section{1}{Introduction}
Consider a canonical system of differential equations in a
neighbourhood of an equilibrium, with Hamiltonian
$$
H(x,y) = H_0(x,y) + H_1(x,y) + \ldots \ ,\quad
(x,y)\in\complessi^{2n}\ , 
\formula{ham.1}
$$
where the unperturbed quadratic part of the Hamiltonian is 
$$
H_0(x,y) = \sum_{j=1}^n \lambda_j x_jy_j\ ,\quad 
(\lambda_1,\ldots,\lambda_n)\in\complessi^n\ ,
\formula{ham.2}
$$
and $H_s(x,y)$ for $s\ge 1$ is a homogeneous polynomial of degree
$s+2$.  The form~\frmref{ham.2} is a typical one for the quadratic
part of a Hamiltonian system in the neighbourhood of an equilibrium,
as under quite general conditions the system may be given that form
via a (complex) linear canonical transformation (see,
e.g.,~\dbiref{Lyapounov-1892} or~\dbiref{Siegel-1971}, \S$\,$15).

The Hamiltonian is assumed to be analytic in some neighborhood of the
origin of $\complessi^{2n}$.  Moreover $\lambda_1$ will be assumed to
satisfy at least the first of the following  non--resonance conditions: 

\item{(i)}{\corsivo First Melnikov's condition:} 
$$
\lambda_{\nu}- k\lambda_1 \ne 0
\quad{\rm for}\ k\in\interi \ {\rm and}\ \nu=1,\ldots,n\>.
\formula{meln.1}
$$

\item{(ii)}{\corsivo Second Melnikov's condition:}
$$
\lambda_{\nu}\pm\lambda_{\nu'}-k\lambda_1 \ne 0
\quad{\rm for}\  k\in\interi \ {\rm and}\ \nu,\nu'=1,\ldots,n\>,
\formula{meln.2}
$$
the case $\nu'=\nu$ being included.

\noindent
In~\dbiref{Lyapounov-1892} Lyapounov proved that if
$\lambda_1=i\omega_1$ is pure imaginary and the non resonance
condition~(i) above is satisfied then there exists a two parameter
family of solutions of the form
$$
x_{j}=\phi_j(\xi_1,\eta_1)\ ,\quad
y_{j}=\psi_j(\xi_1,\eta_1)
\formula{1.1}
$$
written as convergent power series in the arguments 
$$
\xi_1 =
 \oxi_{1}e^{ita_1(\ozeta)}\ ,\quad 
  \eta_1 = 
   \oeta_1e^{-ita_1(\ozeta)}\ ,
\formula{1.2}
$$
where $a_1(\ozeta)=\lambda_1+\ldots\,$
is a convergent power series in $\ozeta_1=\oxi_1\oeta_1$.  In the case
$n=1$ this actually describes all solutions of the system.  A proof in
case all $\lambda$'s are pure imaginary is reported
in~\dbiref{Siegel-1971}.

The proof of the theorem is worked out by the authors quoted above by
expanding the solution in the form~\frmref{1.1} and proceeding by
comparison of coefficients.  From a formal
viewpoint the statement above looks equivalent to the existence of a
canonical transformation that gives the system~\frmref{ham.1} a
suitable normal form, making the Hamiltonian to depend at least
quadratically on $x_2,\ldots,x_n,y_2,\ldots,y_n\,$.  A
formal construction giving such a normal form can be easily produced.
However, proving the convergence of the normalization procedure seems
to be more difficult.  The aim of this paper is precisely to produce a
proof of convergence of the transformation to normal form.

I will actually give two different statements that can be proved with
the same method.  The first one is

\theorem{nrm.1}{With the nonresonance hypothesis~(i) above (first
Melnikov's condition) on
$\lambda_1,\ldots,\lambda_n$, there exists a canonical, near the
identity transformation in the form of a power series convergent in a
neighbourhood of the origin, which gives the
Hamiltonian~\frmref{ham.1} the normal form
$$
H(x,y) = H_0(x,y) + \Gamma(x_1 y_1)+ F(x,y) \ ,\quad
\formula{lyap.1}
$$
where $H_{0}(x,y)$ as in~\frmref{ham.1}, $\Gamma(x_1 y_1)$ depends only on
the product $x_1 y_1$, and $F(x,y)$ is at least quadratic in
$x_2,\ldots,x_n,y_2,\ldots,y_n$}\endclaim

The existence of the Lyapounov orbits for $\lambda_1$ pure imaginary
is evident from the normal form: just put initially
$x_2=\ldots=x_n=y_2=\ldots=y_n=0$, which defines a local invariant two
dimensional manifold on which the dynamics is generated by the
Hamiltonian $\lambda_1 x_1y_1 +\Gamma(x_1 y_1)$.  The advantage of the
normal form is that it allows also to investigate the dynamics in the
neighbourhood of the orbits so found.  To this end the following
statement may be even more useful.

\theorem{nrm.1.a}{With the nonresonance hypotheses~(i) and~(ii) above
(first and second Melnikov's conditions) on
$\lambda_1,\ldots,\lambda_n$, there exists a canonical, near the
identity transformation in the form of a power series convergent in a
neighbourhood of the origin, which gives the
Hamiltonian~\frmref{ham.1} the normal form
$$
H(x,y) = H_0(x,y) + \Gamma(x_1 y_1,\ldots, x_n y_n)+ F(x,y) \ ,\quad
\formula{lyap.1}
$$
where $H_{0}(x,y)$ as in~\frmref{ham.1}, $\Gamma(x,y)$ contains only
monomials of the form $x_1^jy_1^jx_{\nu}y_{\nu}$ with a positive
integer $j$ and with $\nu=2,\ldots,n$, and $F(x,y)$ is at least cubic
in $x_2,\ldots,x_n,y_2,\ldots,y_n$}\endclaim

This requires a stronger non--resonance condition.  However this
normal form may be more convenient if one is interested in the
stability of a Lyapounov orbit.  Indeed, let all $\lambda$ be pure
imaginary, say $\lambda_j=i\omega_j$, and write the Hamiltonian
restricted to the invariant manifold $x_2=\ldots=x_n=y_2=\ldots=y_n=0$
in action--angle variables by transforming
$x_1=\sqrt{p}\,e^{iq}\,,\>y_1=-i\sqrt{p}\,e^{iq}$.  Thus one gets the
Hamiltonian $\omega_1p +\Gamma(p)$, which represents a non linear
oscillator, with orbits written as $p(t)=p_*\,,\>q(t)=q(0) \Omega(p_*)
t$, where $\Omega(p_*) = \omega_1+O(p_*)$ is a fixed frequency.  By a
translation $p'=p-p_*$ the Hamiltonian may be reexpanded (omitting
primes) as
$$
H(q,p,x,y) = \Omega p + \sum_{j=2}^{n} \lambda_jx_jy_j + H_1 + H_2 +\ldots
$$
where $H_s$ is a homogeneous polynomial of degree $s+2$ in
$p^{1/2},x_2,\ldots,y_n$ with coefficients periodically depending on
$q$.  The dynamics of the latter Hamiltonian may be investigated with
known methods from perturbation theory.  The advantage with respect to
the normal form of theorem~\thrref{nrm.1} is that the quadratic part
of the Hamiltonian is independent of the angle $q$.

The proof is based on a previous work by the
author~\dbiref{Giorgilli-2001} concerning the construction of the
normal form in a case investigated by Cherry~\dbiref{Cherry-1926} and
Moser~\dbiref{Moser-1958}.  It must be stressed that this problem does
not involve small divisors.  Rather, the possible source of divergence
is due to the use of Cauchy's estimates for the derivatives required
by the normalization algorithm.  The global effect of accumulation of
derivatives is controlled with a technique introduced by the author
and U. Locatelli in order to achieve a proof of KAM theorem using
classical expansions in a perturbation parameter
(see~\dbiref{Giorgilli-97.3},\dbiref{Giorgilli-97.4},\dbiref{Giorgilli-98a},\dbiref{Giorgilli-99}).

\section{2}{Formal algorithm}
Reducing the Hamiltonian to a normal form is a quite general problem
which may be solved in a number of different ways.  Moreover, the
concept of ``normal form'' may assume a quite general meaning,
depending on what one is looking for.  Here I state the algorithm in a
general form, using the method of composition of Lie series.

\subsection{2.1}{The algorithm for the normal form}
Write the Hamiltonian after $r$ normalization steps as
$$
H^{(r)}(x,y) = H_0(x,y) + Z_1(x,y)+\ldots+Z_r(x,y) 
              +\sum_{s>r} H^{(r)}_s(x,y)\ ,
\formula{alg.1}
$$
where $Z_1(x,y),\ldots,Z_r(x,y)$ are in normal form, whatever it
means, and are homogeneous polynomials of degree $3,\ldots,r+2$.  For
$r=0$ the Hamiltonian~\frmref{ham.1} is considered to be already in
the wanted form, with no functions $Z$.

Assume that the Hamiltonian has been given a normal
form~\frmref{alg.1} up to order $r-1$, so that $H^{(r-1)}$ is known.
The generating function $\chi_r$ and the normal form $Z_r$ are
determined by solving the equation
$$
\lie{H_0}\chi_r + Z_r = H^{(r-1)}_r\ .
\formula{alg.2}
$$
where the common notation $\lie{\phi}\cdot:=\poisson{\cdot}{\phi}$ has
been used.  The solution of this equation depends on what is meant by
``normal form''.  At a formal level, any form that allows to
solve the equation above for $Z_r$ and $\chi_r$ is acceptable.  Assume
for a moment that a method of solution has been found.  Then the
transformed Hamiltonian is expanded as 
$$
\vcenter{\openup1\jot\halign{
\hfil$\displaystyle{#}$
&$\displaystyle{#}$\hfil
&$\displaystyle{#}$\hfil
\cr
%1
H^{(r)}_{sr+m} 
&= \frac{1}{s!}\lie{\chi_r}^{s} Z_m 
   +\sum_{p=0}^{s-1} \frac{1}{p!}\lie{\chi_r}^p H^{(r-1)}_{(s-p)r+m}
     \quad
%\cr
%
%1.a
&{\rm for}\ r\ge 2 \>,\, s\ge 1\ {\rm and}\ 1\le m\lt r\ ,\quad
\cr
%2
H^{(r)}_{sr} 
&= \frac{1}{(s-1)!}\lie{\chi_r}^{s-1} 
    \left(\frac{1}{s} Z_r + \frac{s-1}{s} H^{(r-1)}_{r}
     \right)
     +\sum_{p=0}^{s-2} \frac{1}{p!}\lie{\chi_r}^p H^{(r-1)}_{(s-p)r}
       \hidewidth
\cr
%2.a
&&{\rm for}\ r\ge 1 \ {\rm and}\ s\ge 2\ .\quad
\cr
}}
\formula{alg.3}
$$
The justification of the algorithm requires only some straightforward
calculation, and is deferred to appendix~\appref{lyap.app.a}.

Thus, the problem is how to solve the equation~\frmref{alg.2} for the
generating function and the normal form.  Let me make some general
considerations.

Let $\Pscr_s$ denote the linear space of homogeneous
polynomials of degree $s$ in the complex variables $x,y$. Let also
$\Pscr=\bigcup_{s\ge 0} \Pscr_s$, so that a formal power series is an
element of $\Pscr$.  A basis in $\Pscr$ is given by the monomials
$x^jy^k:=x_1^{j_1}\cdot\ldots\cdot x_n^{j_n}y_1^{k_1}\cdot\ldots\cdot
y_n^{k_n}$, where $j,\,k$ are integer vectors with non--negative
components.
The linear operator $\lie{H_0}$ maps every space $\Pscr_s$ into
itself.  If, due to the choice of the coordinates, the unperturbed
Hamiltonian $H_0$ has the form~\frmref{ham.2} then the operator
$\lie{H_0}$ is diagonal, since
$$
\lie{H_0} x^j y^k = \langle j-k,\lambda\rangle x^j y^k\ .
$$ 
The kernel and the range of $\lie{H_0}$ are defined as usual, namely
$\Nscr=\lie{H_0}^{-1}(0)$, the inverse image of the null vector in
$\Pscr$, and $\Rscr=\lie{H_0}(\Pscr)$. Both $\Nscr$ and $\Rscr$ are
actually subspaces of the same space $\Pscr$, and it turns out that
they are complementary subspaces, i.e., $\Nscr \cap \Rscr = \{0\}$,
the null vector, and $\Nscr \oplus \Rscr = \Pscr$.  A consequence of
the properties above is that $\lie{H_0}$ restricted to the subspace
$\Rscr$ is uniquely inverted, i.e., the equation $\lie{H_0}\chi=\psi$
with $\psi\in\Rscr$ admits an unique solution $\chi$ satisfying the
condition $\chi\in\Rscr$.  That unique solution will be written as
$\chi=\lie{H_0}^{-1}\psi$, i.e., $\lie{H_0}^{-1}$ is defined as the
inverse of $\lie{H_0}$ restricted to $\Rscr$.  It's easy to identify
the subspaces $\Nscr$ and $\Rscr$ using the coordinates.  Thanks to
the diagonal form of $\lie{H_0}$ one has
$$
\eqalign{
\Nscr
&=\span\left\{x^jy^k\>:\> \langle j-k,\lambda\rangle =0\right\}\ ,
\cr
\Rscr
&=\span\left\{x^jy^k\>:\> \langle j-k,\lambda\rangle \ne 0\right\}\ .
\cr
}
\formula{nucleo-range}
$$
Given $\psi\in\Rscr$ and writing $\psi=\sum_{j,k}\psi_{j,k}x^jy^k$,
with $\psi_{j,k}=0$ for $x^jy^k\in\Nscr$, one has
$$
\lie{H_0}^{-1}\psi = \sum_{j,k} \frac{\psi_{j,k}}{\langle
j-k,\lambda\rangle} x^jy^k\ .
\formula{gen_soluzione}
$$

In view of the general considerations above we can conclude that the
choice of a normal form is subjected to the constraint that in
equation~\frmref{alg.2} we have $H^{(r-1)}_r-Z_r\in\Rscr$.  The
simplest choice is to ask also $Z_r\in\Nscr$, i.e., to set $Z_r$ to
be the projection of $H^{(r-1)}_r$ on the subspace $\Nscr$.  This is
known indeed as Birkhoff's normal form.

\subsection{2.2}{Normal form for Lyapounov's orbits}
I come now to show that the construction of the normal form of
theorem~\thrref{nrm.1} is formally consistent.  Consider the disjoint
subsets of $\interi^n$
$$
\eqalign{
\Kscr^{\sharp}
&= \left\{k\in\interi^n\>:\>k_2=\ldots=k_n=0\right\}\ ,
\cr
\Kscr^{\natural}
&=\left\{k\in\interi^n\>:\> |k_2|+\ldots+|k_n|=1\right\}\ ,
\cr
\Kscr^{\flat}
&=\left\{k\in\interi^n\>:\>|k_2|+\ldots+|k_n|>1\right\}\ .
\cr
}
\formula{indici}
$$
One has
$\interi^n=\Kscr^{\sharp}\cup\Kscr^{\natural}\cup\Kscr^{\flat}$, of
course.  Considering only integer vectors $j,\,k$ with non--negative
components, introduce the subspaces of $\Pscr$
$$
\eqalign{
\Pscr^{\sharp}
&=\span\left\{x^jy^k\>:\> j+k\in\Kscr^{\sharp}\right\}
\cr
\Pscr^{\natural}
&=\span\left\{x^jy^k\>:\> j+k\in\Kscr^{\natural}\right\}
\cr
\Pscr^{\flat}
&=\span\left\{x^jy^k\>:\> j+k\in\Kscr^{\flat}\right\}
\cr
}
\formula{sottospazi}
$$
These subspaces are clearly disjoint, and moreover one has
$\Pscr=\Pscr^{\sharp}\oplus\Pscr^{\natural}\oplus\Pscr^{\flat}$.
Finally, let $\Nscr^{\sharp}=\Nscr\cap\Pscr^{\sharp}$ and
$\Rscr^{\sharp}=\Rscr\cap\Pscr^{\sharp}$, and define the subspaces
$\Zscr$ and $\Wscr$ of $\Pscr$ as
$$
\Zscr=\Nscr^{\sharp}\oplus\Pscr^{\flat}\ ,\quad
\Wscr=\Rscr^{\sharp}\oplus\Pscr^{\natural}\ .
\formula{splitting}
$$
It is an easy matter to check that $\Zscr\cap\Wscr=\{0\}$ and
$\Zscr\oplus\Wscr=\Pscr$.  The construction of bases for $\Zscr$ and
$\Wscr$ is quite straightforward: a monomial $x^jy^k$ belongs to
$\Zscr$ in either case ($j+k\in\Kscr^{\sharp}$ and $\langle
j-k,\lambda\rangle=0$) or ($j+k\in\Kscr^{\flat}$); else it belongs to
$\Wscr$.  The hypothesis~(i) on $\lambda$ (first Melnikov's condition)
formulated at the beginning of the introduction means that the
non-resonance condition
$$
\langle k,\lambda\rangle \ne 0\quad
\hbox{for $0\ne k\in \Kscr^{\sharp}\cup\Kscr^{\natural}$}
\formula{nonrisonanza}
$$
is satisfied. This implies $\Wscr\subset\Rscr$, so that for every
$\psi\in\Wscr$ the unique solution $\chi=\lie{H_0}^{-1}\psi$,
$\chi\in\Wscr$ of the equation $\lie{H_0}\chi=\psi$ exists.  With this
setting, the equation
$$
\lie{H_0} \chi + Z = \Psi\ ,
\formula{gen.1}
$$
with $\Psi$ known, admits a straightforward solution.  Split
$\Psi=\Psi_{\Zscr}+\Psi_{\Wscr}$ with $\Psi_{\Zscr}\in\Zscr$ and
$\Psi_{\Wscr}\in\Wscr$; such a decomposition exists and is unique,
because $\Zscr$ and $\Wscr$ are complementary subspaces.  Then set
$Z=\Psi_{\Zscr}$, and determine $\chi=\lie{H_0}^{-1}\Psi_{\Wscr}$
according to~\frmref{gen_soluzione}.

\subsection{2.3}{Adding the second Melnikov's condition}
With minor changes one can also prove that the normal form of
theorem~\thrref{nrm.1.a} can be constructed.  Let
$$
\eqalign{
\Kscr^{\sharp}
&= \left\{k\in\interi^n\>:\>k_2=\ldots=k_n=0\right\}\ ,
\cr
\Kscr^{\natural}
&=\left\{k\in\interi^n\>:\> |k_2|+\ldots+|k_n|=1,2\right\}\ ,
\cr
\Kscr^{\flat}
&=\left\{k\in\interi^n\>:\>|k_2|+\ldots+|k_n|>2\right\}\ .
\cr
}
\formula{indici.a}
$$
One has again
$\interi^n=\Kscr^{\sharp}\cup\Kscr^{\natural}\cup\Kscr^{\flat}$, of
course.  The subspaces of $\Pscr$ are defined again as
in~\frmref{sottospazi}, although they turn out to be different in view
of the differences in the sets $\Kscr$.  Finally, let
$\Nscr^{\sharp}=\Nscr\cap\Pscr^{\sharp}$,
$\Rscr^{\sharp}=\Rscr\cap\Pscr^{\sharp}$,
$\Nscr^{\natural}=\Nscr\cap\Pscr^{\natural}$ and
$\Rscr^{\natural}=\Rscr\cap\Pscr^{\natural}$, and define the subspaces
$\Zscr$ and $\Wscr$ of $\Pscr$ as
$$
\Zscr=\Nscr^{\sharp}\oplus\Nscr^{\natural}\oplus\Pscr^{\flat}\ ,\quad
\Wscr=\Rscr^{\sharp}\oplus\Rscr^{\natural}\ .
\formula{splitting.a}
$$
The difference with respect to the previous case is just that now
$\Nscr^{\natural}$ is not empty, because it contains all monomials of
the form $(x_1y_1)^j\times x_{\nu}y_{\nu}$ with $\nu=2,\ldots,n$ and
positive $j$.  This forces the change in the definition of the
subspaces $\Zscr$ and $\Wscr$.  However, the properties
$\Zscr\cap\Wscr=\{0\}$ and $\Zscr\oplus\Wscr=\Pscr$ remain true.
Furthermore, in view of the second Melnikov's condition, also the
property $\Wscr\subset\Rscr$ holds true, so that for every
$\psi\in\Wscr$ the unique solution $\chi=\lie{H_0}^{-1}\psi$,
$\chi\in\Wscr$ of the equation $\lie{H_0}\chi=\psi$ exists.

\section{3}{Quantitative estimates}
Pick a real vector $R\in\reali^n$ with positive components. and
consider the domain
$$
\Delta_R = \left\{(x,y)\in\complessi^n
 \>:\> 
  |x_j|\le R_j\,,\>|y_j|\le R_j\>,\,1\le j\le n
   \right\}\ ,
\formula{disk}
$$
namely a polydisk which is the product of disks of radii
$R_1,\ldots,R_n$ in the planes of the complex coordinates
$(x_1,\ldots,x_n)$  and $(y_1,\ldots,y_n)$, respectively. 
Let also
$$
\Lambda = {\min_{1\le j\le n} R_j}\ .
\formula{lambda.disk}
$$

The norm $\ordnorma{f}_{R}$ in the polydisk $\Delta_R$ is defined as
$$
\ordnorma{f}_R = \sum_{|j+k|=r} |f_{j,k}| R^{j+k}\ .
\formula{norm.pol}
$$
A family of polydisks $\Delta_{\delta R}$ of radii $\delta R$, with
$0\lt\delta\le 1$ will be considered below.  With a minor abuse the
simplified notation $\|\cdot\|_{\delta}$ in place of
$\|\cdot\|_{\delta R}$ will be used.

The main result of this section is 

\lemma{generatrici}{Let the Hamiltonian $H^{(0)}$ satisfy
$\ordnorma{H^{(0)}_s}_1\le h^{s-1}E$ for $s\ge 1$, with some constants
$h\ge 0$ and $E\gt 0$. Let $0\lt d\lt 1/2$. Then there exist positive
constants $\beta$ and $G$ depending on $E,\,h,\,\Lambda,\,d$ and on
$\lambda_1,\ldots,\lambda_n$ such that
$$
\ordnorma{\chi_r}_{1-d}\le \beta^{r-1}G
 \quad \hbox{for all $r\ge 1$}\ .
$$
}\endclaim

\noindent
The rest of this section is devoted to the proof.  Some technical
calculations are deferred to appendix~\appref{lyap.app.b}.

\subsection{3.2}{An arithmetic lemma}
The following lemma will play a crucial role in the proof of
lemma~\lemref{generatrici}.

\lemma{art.1}{Let $\lambda\in\complessi^n$ be such that
$\lambda_1$ satisfies the non-resonance condition~\frmref{meln.1}.
Then there exists a positive $\gamma$ such that
the inequality
$$
|\langle k,\lambda\rangle | \ge |k|\gamma
$$ 
holds true for all non--zero $k\in\Kscr^{\sharp}\cup\Kscr^{\natural}$
defined as in~\frmref{indici}.
}\endclaim

\corollary{art.1a}{Let in addition the non-resonance
condition~\frmref{meln.2} be satisfied.  Then the same statement holds
true for all non--zero $k\in\Kscr^{\sharp}\cup\Kscr^{\natural}$
defined as in~\frmref{indici.a}.
}\endclaim

The proof of the corollary is just a trivial modification of the 

\prooftx{of lemma~\lemref{art.1}}
For $k\in\Kscr^{\sharp}$ the claim is obvious, since $|\langle
k,\lambda\rangle |=|k_1\lambda_1|$. So, let $k\in\Kscr^{\natural}$.  Set
$\theta=\max\bigl(|\lambda_2|,\ldots,|\lambda_n|\bigr)$ (for the
corollary maximize also over $|2\lambda_{\nu}|$ and
$|\lambda_{\nu}\pm\lambda_{\nu'}|$ with $\nu'\ne\nu$).  Pick an
integer $N\ge 1+2\theta$ and set
$$
\delta=\min_{{k\in\Kscr^{\natural}}\atop{|k|\le N}} \bigl|\langle
k,\lambda\rangle\bigr| 
\ ,\quad \gamma = \min \left(\frac{\delta}{N},\frac{|\lambda_1|}{2}\right)
\ ;
$$
in view of the non-resonance condition~\frmref{nonrisonanza} one has
$\delta\gt 0$.  Then the claim of the lemma holds true with the given
value of $\gamma$.  Indeed, let $k\in\Kscr^{\natural}$, so that
$|k_1|=|k|-1$.  If $|k|\le N$ then $\bigl|\langle
k,\lambda\rangle \bigr|\ge \delta\ge N\gamma \ge |k|\gamma$.  If
$|k|\gt N$ use $\theta\le (N-1)\delta/2$, which follows from the
choice of $N$, and evaluate
$$
\displaylines{
\qquad
 \bigl|\langle k,\lambda\rangle \bigr|
  \ge \bigl|k_1\lambda_1\bigr| - \theta
   \ge \bigl(|k|-1 \bigr)\delta - \frac{(N-1)}{2}\,\delta
\hfill\cr\hfill
\ge \frac{|k|-1}{2}\,\delta + \frac{N}{2}\,\delta -
     \frac{(N-1)}{2}\,\delta
 = |k|\frac{\delta}{2} \ge |k|\gamma\ .
\qquad
}
$$
\endproof

\subsection{3.1}{Generalized Cauchy estimates}
Here I refer to the more restrictive hypotheses of
theorem~\thrref{nrm.1.a}, and in particular to the spaces $\Pscr$
defined as in sect.~\sbsref{2.3}.  However, the same arguments with
very little simplifications apply also to the setting of
sect.~\sbsref{2.2}, which applies to theorem~\thrref{nrm.1}.  I will
insert short comments in parentheses concerning the latter case, where
appropriate.

The estimates in this section strongly depend on a suitable
splitting of all functions over the subspaces
$\Pscr^{\sharp},\,\Pscr^{\natural}$ and $\Pscr^{\flat}$.  At a formal
level, it is useful to keep in mind the following table concerning the
Poisson bracket:
$$
\vcenter{\tabskip=0pt
\def\tablerule{\noalign{\hrule}}
\halign{
 \hbox to 3 em{\hfil$\displaystyle{#}\hfil$}\vrule\vrule
&\hbox to 6.5 em{\hfil$\displaystyle{#}\hfil$}\bigg\vert
&\hbox to 6.5 em{\hfil$\displaystyle{#}\hfil$}\bigg\vert
&\hbox to 6.5 em{\hfil$\displaystyle{#}\hfil$}\vrule
\cr
\{{\cdot},{\cdot}\}
 & \Pscr^{\sharp}
 & \Pscr^{\natural}
 & \Pscr^{\flat}
\cr
\tablerule
\tablerule
\Pscr^{\sharp}
 & \Pscr^{\sharp}
 & \Pscr^{\natural}
 & \Pscr^{\flat}
\cr
\tablerule
\Pscr^{\natural}
 & \Pscr^{\natural}
 & \Pscr^{\sharp} \oplus \Pscr^{\natural} \oplus \Pscr^{\flat}
 & \Pscr^{\natural} \oplus \Pscr^{\flat}
\cr
\tablerule
\Pscr^{\flat}
 & \Pscr^{\flat}
 & \Pscr^{\natural} \oplus \Pscr^{\flat}
 & \Pscr^{\flat}
\cr
\tablerule
}}
\formula{tabella_poisson}
$$
(For the subspaces defined as in sect.\sbsref{2.2} just remove
$\Pscr^{\natural}$ from the central case corresponding to the Poisson
bracket between functions in $\Pscr^{\natural}$.)  In view of the
transformation formul\ae~\frmref{alg.3} the situation to be considered
is the following.  A generating function $\chi\in\Wscr\cap\Pscr_r$
with some $r\ge 1$ is given in the form $\chi=\lie{H_0}^{-1}\psi$,
with known $\psi\in\Wscr\cap\Pscr_r$.  Since
$\Wscr=\Rscr^{\sharp}\oplus\Rscr^{\natural}$ one has
$\chi=\chi^{\sharp}+\chi^{\natural}$, with an obvious meaning of the
notation.  The operator $\lie{\chi}$ may be applied either to a
generic function $f=f^{\sharp}+ f^{\natural} +f^{\flat}\in\Pscr_s$
with $s\ge r$ or to a function in normal form
$Z=Z^{\sharp}+Z^{\natural}+Z^{\flat}\in\Zscr\cap\Pscr_m$ with $0\lt
m\lt r$ (in theorem~\thrref{nrm.1} $Z^{\natural}=0$).  In particular
one has 
$$
Z^{\sharp}=\sum_{j>1} z_{j}x_1^{j}y_1^{j}\ ,\quad
Z^{\natural}=\sum_{j\gt 0,\,2\le\nu\le n}
z_{j,\nu}x_1^{j}y_1^{j}x_{\nu}y_{\nu}\ ,
$$
due to the non--resonance
conditions on $\lambda$.  For some non--negative $\delta',\delta'',\delta$ satisfying
$0\le\max(\delta',\delta'')\lt \delta\le 1/2$ the norms $\norma{\psi}_{1-\delta'}\,$,
$\norma{f}_{1-\delta''}$ and $\norma{Z}_{1-\delta''}$ are assumed to be known,
and one looks for an estimate of the Lie derivative in a domain
$\Delta_{(1-\delta)R}$.  The following estimates will be used in the rest
of the paper.

\item{(i)}The generating function $\chi$ is estimated by
$$
\norma{\chi}_{1-\delta'} \le \frac{1}{\gamma} \norma{\psi}_{1-\delta'}\ ,
\formula{stimagen}
$$
with $\gamma$ as in lemma~\lemref{art.1}.

\item{(ii)}The general estimate for the Lie derivative of a generic
function $f$ is
$$
\bignorma{\lie{\chi}f}_{1-\delta} 
 \le \frac{1}{(\delta-\delta')(\delta-\delta'')\Lambda^2} 
  \norma{\chi}_{1-\delta'}
   \norma{f}_{1-\delta''}
\formula{cauchy.f}
$$
with $\Lambda$ as in~\frmref{lambda.disk}.  Denoting by
$\bigl(\lie{\chi^{\natural}}f^{\flat}\bigr)^{\natural}$ the projection
of $\lie{\chi^{\natural}}f^{\flat}$ over $\Pscr^{\natural}$ one has
$$
\bignorma{\bigl(\lie{\chi^{\natural}}f^{\flat}\bigr)^{\natural}}_{1-\delta}
 \le \frac{4}{(\delta-\delta'')\Lambda^2} 
  \norma{\chi}_{1-\delta'}
   \norma{f}_{1-\delta''}\ .
\formula{cauchy.f-beq}
$$

\item{(iii)}For a function $Z$ in normal form one has
$$
\bignorma{\lie{\chi}(Z^{\sharp}+Z^{\natural})}_{1-\delta}
 \le \frac{1}{(\delta-\delta'')\gamma\Lambda^2} 
  \norma{\psi}_{1-\delta'}
   \norma{Z}_{1-\delta''}\ .
\formula{cauchy.Z-dies}
$$

\noindent
I recall the reader's attention on the missing denominator $\delta-\delta'$
in~\frmref{cauchy.f-beq} and~\frmref{cauchy.Z-dies}.  This is crucial
for the convergence proof.  For, working out the convergence proof
requires a quite accurate control of the accumulation of the divisors
$\delta-\delta',\, \delta-\delta''$ that appear in the generalized Cauchy estimates for
derivatives.  The scheme in the next section is specially devised in
order to allow such a control.

The proof of~\frmref{stimagen} is a straightforward consequence of the
definition of the norm and of~\frmref{gen_soluzione}.  For, the
denominators are uniformly estimated from below by $\gamma$, in view of
lemma~\lemref{art.1}. 

The proof of the estimates~\frmref{cauchy.f}, \frmref{cauchy.f-beq}
and~\frmref{cauchy.Z-dies} is a purely technical matter, and is
deferred to appendix~\appref{lyap.app.b}.

\subsection{3.3}{Recursive estimates}
The aim of this section is to obtain estimates for the norms of the
generating functions and of the transformed Hamiltonians, at every step
of the normalization procedure.

Consider a sequence of boxed domains $\Delta_{(1-\delta_r)R}$, where
$\{\delta_r\}_{r\ge 1}$ is a monotonically increasing sequence of
positive numbers converging to some $d\lt 1/2$.  Let also
$\delta_0=0$, and $d_r=\delta_r-\delta_{r-1}$ for $r\ge 1$, so that
$d_r\lt 1$ for all positive $r$.  The purpose is to look for estimates
of the norms of the generating function $\chi_r$ and of the normal
form $Z_r$ in the polydisk $\Delta_{(1-\delta_{r-1})R}$, and of the
functions $H^{(r)}_s$ in the domain $\Delta_{(1-\delta_r)R}$.

Let $\Jscr_{r,s}$ for $1\lt r\lt s$ be the set of integer arrays
defined as
\formdef{jset}
$$
\displaylines{
\frmref{jset}\qquad
\Jscr_{r,s} = 
 \Bigl\{J=\{j_1,\ldots j_k\}\>:\> 
         j_m\in\{1,\ldots,r\}\,,\>
          1\le k\le 2(s-1)\,,\>
 \Bigr.\hfill\cr\hfill\Bigl.
           \sum_{m=1}^{k} \lgd j_m\le 2(s-1-\lgd s)
 \Bigr\}\ .
\qquad\cr
}
$$
Let also $\Jscr_{0,s}=\emptyset$ for $s\ge 1$.  Recalling that
$\{d_r\}_{r\ge 1}$ is a sequence of positive numbers not exceeding $1$
define the sequence $\{T_{r,s}\}_{0\le r\lt s}$ as
$$
T_{0,s}=1\ ,\quad
T_{r,s}= \max_{J\in\Jscr_{r,s}} \prod_{j\in J} d_j^{-1}\ .
\formula{T}
$$
The following properties will be used below: for $0\le r\le r'\lt s$
one has
\formdef{T-prop.1}
\formdef{T-prop.2}
$$
\eqalignno{
T_{r,s} &\le T_{r',s}\ , &\frmref{T-prop.1}\cr
\frac{1}{d_r^2} T_{r-1,r} T_{r',s}
&\le T_{r',r+s}\ . &\frmref{T-prop.2}\cr
}
$$
Checking~\frmref{T-prop.1} is easy: for $r=0$ use $d_l\le
1$ for $l\ge 1$; for $r\gt 0$ use the inclusion relation
$\Jscr_{r,s}\subset\Jscr_{r',s}$ for $r\lt r'$.  In order to
prove~\frmref{T-prop.2} remark that by definition one has
$$
\eqalign{
\frac{1}{d_r^2} T_{r-1,r} T_{r',s}
&= \frac{1}{d_r^2} 
    \max_{J\in\Jscr_{r-1,r}} \prod_{j\in J} d_j^{-1}
     \max_{J'\in\Jscr_{r',s}} \prod_{j'\in J'} d_{j'}^{-1}
\cr
& = \max_{J\in\Jscr_{r-1,r}} \max_{J'\in\Jscr_{r',s}}
     \prod_{j\in\{r,r\}\cup J \cup J'} d_j^{-1}\ .
\cr
}
$$
It is enough to prove that $\{r,r\}\cup J \cup J' =:\tilde J \in
\Jscr_{r',r+s}$.  First check that
$$
\#\bigl(\tilde J) =
 2 + \#(J) +\#(J')
  \le 2 + 2(r-1) + 2(s-1) = 2(r+s-1)\ .
$$
On the other hand, since $1\le j \le r-1$ for all $j\in J$ and $1\le j'
\le r'$ for all $j'\in J'$, one also has $1 \le \tilde j \le r'$ for
all $\tilde j \in \tilde J$.  Finally, evaluate
$$
\eqalign{
\sum_{\tilde j\in\tilde J} \lgd \tilde j
&= 2 \lgd r + \sum_{j\in J} \lgd j +
        \sum_{j'\in J'} \lgd j'
\cr
& \le 2 \lgd r + 2(r-1-\lgd r) + 2(s-1-\lgd s)
\cr
& \le 2\bigl[r+s-1-(1+\lgd s)\bigr]
\cr
& \le 2\bigl[r+s-1-\lgd(r+s)\bigr]\ ,
}
$$
where the elementary inequality $1+\lgd s=\lgd 2 + \lgd s = \lgd (2s)
\gt \lgd (r+s)$ has been used (recall that $r\le r'\lt s$).  Hence,
$\tilde J\in\Jscr_{r',r+s}$, as claimed.

I shall also use the numerical sequence $\{\mu_{r,s}\}_{r\ge 0,s\ge
0}$ defined as
$$
\eqalign{
\mu_{0,0}&=0\>,\quad\mu_{0,s}=1\quad\hbox{for $s\gt 0$}\ ,
\cr
\mu_{r,s} &= \sum_{0\le rp\lt s} \mu_{r-1,r}^p \mu_{r-1,s-rp}
 \quad
  \hbox{ for $r\gt 0$ and $s\ge 0$}\ .
}
\formula{museq}
$$

The recursive estimates are collected in 

\lemma{nrm.1}{Let the Hamiltonian $H^{(0)}$ satisfy
$\ordnorma{H^{(0)}_s}_1\le h^{s-1}E$ for some constants $h\ge 0$ and
$E\gt 0$.  Let $d_0=1$ and $\{d_r\}_{r\ge 1}$ be an arbitrary sequence of positive
numbers satisfying $\sum_{r\ge 1}d_r=d$ with $d\lt 1$. Let also $\delta_0=0$ and 
$\delta_r=d_1+\ldots+d_r$.  Then for $s\gt r\ge 1$ the following
estimates hold true:
\formdef{nrm.1}
\formdef{nrm.2}
\formdef{nrm.2a}
\formdef{nrm.3}
\formdef{nrm.3a}
$$
\eqalignno{
\ordnorma{\chi_r}_{1-\delta_{r-1}} 
&\le 
 \mu_{r-1,r} T_{r-1,r} C^{r-1} \frac{E}{\gamma} \ ,
&\frmref{nrm.1}
\cr
\ordnorma{Z_r}_{1-\delta_{r-1}}
&\le 
 \mu_{r-1,r} T_{r-1,r} C^{r-1} \frac{E}{d_{r-1}}  \ ,
&\frmref{nrm.2}
\cr
\ordnorma{Z_r^{\sharp}+Z_r^{\natural}}_{1-\delta_{r-1}}
&\le 
 \mu_{r-1,r} T_{r-1,r} C^{r-1} E  \ ,
&\frmref{nrm.2a}
\cr
\ordnorma{H^{(r)}_s}_{1-\delta_{r}} 
&\le 
 \mu_{r,s} T_{r,s} C^{s-1} \frac{E}{d_{r}} \ ,
&\frmref{nrm.3}
\cr
\ordnorma{H^{(r),\sharp}_s+H^{(r),\natural}_s}_{1-\delta_{r}} 
&\le 
 \mu_{r,s} T_{r,s} C^{s-1} E \ ,
&\frmref{nrm.3a}
\cr
}
$$
where
$$
C=h+\frac{4e^2E}{\gamma\Lambda^2}\ ,
\formula{nrm.4}
$$
and $\mu_{r,s}$ and $T_{r,s}$ are the sequences defined by~\frmref{T}
and~\frmref{museq}.}
\endclaim

\noindent
Remark that~\frmref{nrm.1}, \frmref{nrm.2a} and~\frmref{nrm.3a} differ
from~\frmref{nrm.2} and~\frmref{nrm.3}, respectively, only because a
divisor $d_r$ is missing.

\proof
By induction.  For $r=0$ only~\frmref{nrm.3}
and~\frmref{nrm.3a} are meaningful, and hold true in view of
$d_0=\mu_{0,s}=T_{0,s}=1$ and of $h\lt C$.  The induction
consists in first proving that if~\frmref{nrm.3} and~\frmref{nrm.3a}
hold true up to $r-1$ then~\frmref{nrm.1}, \frmref{nrm.2}
and~\frmref{nrm.2a} are true for $r$; next proving that
if~\frmref{nrm.1}, \frmref{nrm.2} and~\frmref{nrm.2a} hold true up to
$r$ then~\frmref{nrm.3} and~\frmref{nrm.3a} are true for $r$.

Let $r\gt 0$ and put $r-1$ in place of $r$ and $r$ in place of $s$
in~\frmref{nrm.3} and~\frmref{nrm.3a}. Recalling that only
$H^{(r-1),\sharp}_r + H^{(r-1),\natural}_r$ is used in order to
determine $\chi_r$ use the definition of the norm, the form of the
solution of eq.~\frmref{alg.2} discussed in sect.~\secref{2}, and the
estimate~\frmref{stimagen}.  This immediately shows
that~\frmref{nrm.1}, \frmref{nrm.2} and~\frmref{nrm.2a} are true for
$r$ provided~\frmref{nrm.3} and~\frmref{nrm.3a} hold true for $r-1$.
Coming to~\frmref{nrm.3} and~\frmref{nrm.3a} and recalling the
recursive definitions~\frmref{alg.3} there are only two kinds of terms
to be estimated, namely $\frac{1}{s!}\lie{\chi_r}^s Z_m$ for $1\le
m\lt r$ and $\frac{1}{p!}\lie{\chi_r}^p H^{(r-1)}_{(s-p)r+m}$ for
$0\le p\le s$ and $0\le m\lt r$.  For, remarking that $Z_r$ and
$H^{(r-1)}_r$ are estimated by exactly the same quantity it is safe to
estimate $\bignorma{H^{(r)}_{sr}}$ by replacing $Z_r$ with
$H^{(r-1)}_r$ in the second of~\frmref{alg.3}.  This is tantamount to
extending the sum in the second of~\frmref{alg.3} to $p=s-1$ and
making it identical with the sum in the first of~\frmref{alg.3}, with
$m=0$.

Denote $\phi_s=\lie{\chi_r}^s Z_m$, where $r\gt 1$, and split $\phi_s
= \phi_s^{\sharp} + \phi_s^{\natural} +
\phi_s^{\flat}$.  I claim  
\formdef{stima.Z}
\formdef{stima.Z.sharp}
$$
\eqalignno{
\bignorma{\phi_s}_{1-\delta_r}
&\le s!\, \left(\frac{e}{d_r \Lambda}\right)^{2(s-1)}
  \norma{\chi_r}_{1-\delta_{r-1}}^{s-1} \frac{D}{d_r}\ ,
&\frmref{stima.Z}
\cr
\bignorma{\phi_s^{\sharp}+\phi_s^{\natural}}_{1-\delta_r}
&\le s!\, \left(\frac{2e}{d_r \Lambda}\right)^{2(s-1)}
  \norma{\chi_r}_{1-\delta_{r-1}}^{s-1} D\ ,
&\frmref{stima.Z.sharp}
\cr
}
$$
for $s\ge 1$, where
$$
D= \mu_{r-1,r} \mu_{m-1,m} T_{r,r+m} C^{r+m-1} E
\formula{D}
$$
The proof proceeds by induction. Let $s=1$. By the general
estimate~\frmref{cauchy.f} one has
$$
\bignorma{\phi_1}_{1-\delta_r} 
 \le \frac{1}{d_r d_m\Lambda^2}
  \bignorma{\chi_r}_{1-\delta_{r-1}}
   \bignorma{Z_m}_{1-\delta_{m-1}}\ .
$$
Using~\frmref{nrm.1} and \frmref{nrm.2} one gets
$$
\bignorma{\phi_1}_{1-\delta_r} 
 \le \frac{1}{d_r}
  \mu_{m-1,m} \mu_{r-1,r} \>
   \frac{1}{d_{m-1}d_m} T_{m-1,m} T_{r-1,r}\> C^{r+m-2} 
    \frac{E^2}{\gamma\Lambda^2}\ ,
$$
so that~\frmref{stima.Z} immediately follows from~\frmref{T-prop.1},
\frmref{T-prop.2} and~\frmref{nrm.4}

Still keeping $s=1$,~\frmref{stima.Z.sharp} is obtained by remarking
that the contributions to $\phi_1^{\sharp}+\phi_1^{\natural}$ come
only from $\lie{\chi_r}(Z_m^{\sharp}+Z_m^{\natural})$ and
$\bigl(\lie{\chi_r^{\natural}}Z_m^{\flat}\bigr)^{\natural}$.
Proceeding as above, from~\frmref{cauchy.Z-dies}
and~\frmref{cauchy.f-beq} one gets
$$
\displaylines{
\qquad\qquad
\bignorma{\phi_1^{\sharp}+\phi_1^{\natural}}_{1-\delta_r}
 \le \frac{1}{\gamma d_m\Lambda^2} 
  \bignorma{H^{(r-1),\sharp}_r + H^{(r-1),\natural}_r}_{1-\delta_{r-1}}
   \bignorma{Z_m}_{1-\delta_{m-1}}
\hfill\cr\hfill
+\frac{4}{d_m\Lambda^2} 
  \bignorma{\chi_r}_{1-\delta_{r-1}}
   \bignorma{Z_m}_{1-\delta_{m-1}}\ .
\qquad\qquad\cr
}
$$
Then~\frmref{stima.Z.sharp} for $s=1$ follows from~\frmref{nrm.2},
\frmref{nrm.3a} for $r-1$ and~\frmref{nrm.4}.  Remark that the divisor
$d_r$ does not appear here.

Let now $s\gt 1$, and assume that~\frmref{stima.Z} be true up to
$s-1$.  Recalling that the divisor $d_r$ due to the generalized Cauchy
estimates is arbitrary, replace $d_r$ with $\frac{s-1}{s}d_r$ in the
estimates~\frmref{stima.Z} and~\frmref{stima.Z.sharp} for
$\phi_{s-1}$, thus getting
$$
\eqalign{
\bignorma{\phi_{s-1}}_{1-\delta_r+d_r/s}
& \le (s-1)!\, \left(\frac{s}{s-1}\right)^{2s-3}
   \left(\frac{e}{d_r \Lambda}\right)^{2(s-2)}
    \norma{\chi_r}_{1-\delta_{r-1}}^{s-2} \frac{D}{d_r}\ ,
\cr
\bignorma{\phi^{\sharp}_{s-1}+\phi^{\natural}_{s-1}}_{1-\delta_r+d_r/s}
& \le (s-1)!\, \left(\frac{s}{s-1}\right)^{2s-4}
   \left(\frac{2e}{d_r \Lambda}\right)^{2(s-2)}
    \norma{\chi_r}_{1-\delta_{r-1}}^{s-2} D\ .
\cr
}
\formula{sta.1}
$$
Consider first the estimate~\frmref{stima.Z.sharp}. Remarking that the
contributions to $\phi_s^{\sharp}+\phi_s^{\natural}$ come only from
$\lie{\chi_r}\bigl(\phi_{s-1}^{\sharp}+\phi_{s-1}^{\natural}\bigr)$
and
$\bigl(\lie{\chi_r^{\natural}}\phi_{s-1}^{\flat}\bigr)^{\natural}$,
use~\frmref{cauchy.f} and~\frmref{cauchy.f-beq} to estimate
$$
\eqalign{
\bignorma{\phi^{\sharp}_{s}+\phi^{\natural}_{s}}_{1-\delta_r}
&\le \bignorma{\lie{\chi_r}\bigl(\phi^{\sharp}_{s-1}+
      \phi^{\natural}_{s-1}\bigr)}_{1-\delta_r} +
       \bignorma{\bigl(\lie{\chi_r^{\natural}}
        \phi_{s-1}^{\flat}\bigr)^{\natural}}_{1-\delta_r}
\cr
&\le \frac{s}{d_r^2\Lambda^2} \bignorma{\chi_r}_{1-\delta_{r-1}}
      \bignorma{\phi^{\sharp}_{s-1}+
       \phi^{\natural}_{s-1}}_{1-\delta_r+d_r/s}
\cr
&\qquad + \frac{2s}{d_r\Lambda^2} 
           \bignorma{\chi_r}_{1-\delta_{r-1}}
            \bignorma{\phi_{s-1}}_{1-\delta_r+d_r/s}
\cr
}
$$
Replacing~\frmref{sta.1} in the latter expression one gets
$$
\bignorma{\phi^{\sharp}_{s}+\phi^{\natural}_{s}}_{1-\delta_r}
 \le \frac{s!}{d_r^2\Lambda^2}
      \left(\frac{s}{s-1}\right)^{2s-3}
       \left(\frac{2e}{d_r\Lambda}\right)^{2(s-2)}
        \bignorma{\chi_r}_{1-\delta_{r-1}}^{s-1} D\ ,
$$
so that~\frmref{stima.Z.sharp} follows from the trivial inequality
$\bigl(\frac{s}{s-1}\bigr)^{s-1}\lt e$.  The estimate~\frmref{stima.Z}
is checked with a similar calculation, just taking into account
that~\frmref{cauchy.f} must be used in order to estimate
$\lie{\chi_r}\phi_{s-1}$.  This produces an extra divisor $d_r$ with
respect to the calculation above.

Finally, replace~\frmref{nrm.1} and~\frmref{D} in~\frmref{stima.Z}
and~\frmref{stima.Z.sharp}. Using also~\frmref{nrm.4},
one gets
$$
\eqalign{
\bignorma{\phi_{s}}_{1-\delta_r}
& \le s!\, \mu_{r-1,r}^{s-1} 
   \left(\frac{1}{d_r^2}\, T_{r-1,r}\right)^{s-1} T_{r,r+m} C^{sr+m-1} 
    \frac{E}{d_r}\ .
\cr
\bignorma{\phi^{\sharp}_{s}+\phi^{\natural}_{s}}_{1-\delta_r}
& \le s!\, \mu_{r-1,r}^{s-1}
   \left(\frac{1}{d_r^2}\, T_{r-1,r}\right)^{s-1} T_{r,r+m} C^{sr+m-1} E\ .
\cr
}
$$
Using $s-1$ times the inequalities~\frmref{T-prop.1}
and~\frmref{T-prop.2} one easily gets
$$
\left(\frac{1}{d_r^2}\, T_{r-1,r}\right)^{s-1} T_{r,r+m}
 \le \left(\frac{1}{d_r^2}\, T_{r-1,r}\right)^{s-2} T_{r,2r+m}
   \le \ldots 
    \le T_{r,sr+m}\ .
$$
Thus one concludes
\formdef{lchi.Z.1}
\formdef{lchi.Z.2}
$$
\eqalignno{
\frac{1}{s!}\,
 \bignorma{\lie{\chi_r}^s Z_m}_{1-\delta_r}
& \le \mu_{r-1,r}^s \mu_{m-1,m} T_{r,sr+m} C^{sr+m-1} 
   \frac{E}{d_r}\ ,
&   \frmref{lchi.Z.1}
\cr
\frac{1}{s!}\,
 \bignorma{\bigl(\lie{\chi_r}^s Z_m\bigr)^{\sharp} + 
            \bigl(\lie{\chi_r}^s Z_m\bigr)^{\natural}}_{1-\delta_r}
&  \le \mu_{r-1,r}^s \mu_{m-1,m} T_{r,sr+m} C^{sr+m-1} E\ .
&   \frmref{lchi.Z.2}
\cr
}
$$

The estimate for $\frac{1}{p!}\lie{\chi_r}^p H^{(r-1)}_{(s-p)r+m}$ is
a minor {\it variazione} of the scheme above.  Only the first step
must be omitted.  E.g., set $\phi_{p} = \lie{\chi_r}^p
H^{(r-1)}_{(s-p)r+m}$ and proceed as follows. Using~\frmref{nrm.3} for
$r-1$ get
$$
\bignorma{\phi_0}_{1-\delta_{r-1}} 
 \le \mu_{r-1,sr+m} T_{r,sr+m} C^{sr+m-1} E\ ;
$$
this starts the induction on $p$.  Then proceed for $p\gt 0$ as
above. The conclusion is
\formdef{lchi.H.1}
\formdef{lchi.H.2}
$$
\vcenter{\openup1\jot\halign to\displaywidth{
{#}\qquad\hfil
&\hfil$\displaystyle{#}$
&$\displaystyle{#}$\hfil
&$\displaystyle{#}$\hfil\tabskip=\centering
\cr
\frmref{lchi.H.1}
&\frac{1}{p!}\,
& \bignorma{\lie{\chi_r}^p H^{(r-1)}_{(s-p)r+m}}_{1-\delta_r}
&  \le \mu_{r-1,r}^p \mu_{r-1,(s-p)r+m} T_{r,sr+m} C^{sr+m-1}
   \frac{E}{d_r}\ , 
\cr
\frmref{lchi.H.2}
&\frac{1}{p!}\,
& \bignorma{\bigl(\lie{\chi_r}^p H^{(r-1)}_{(s-p)r+m}\bigr)^{\sharp} + 
   \bigl(\lie{\chi_r}^p H^{(r-1)}_{(s-p)r+m}\bigr)^{\natural}}_{1-\delta_r}
    \hidewidth 
\cr
& & & \le \mu_{r-1,r}^p \mu_{r-1,(s-p)r+m} T_{r,sr+m} C^{sr+m-1} E\ .
\cr
}}
$$

Collecting~\frmref{lchi.Z.1}, \frmref{lchi.Z.2}, \frmref{lchi.H.1}
and~\frmref{lchi.H.2} and referring to the transformation
formul\ae~\frmref{alg.3} it is now an easy matter to verify
that~\frmref{nrm.3} and~\frmref{nrm.3a} hold true provided the sequence
$\mu_{r,s}$ for $0\lt r\lt s$ is defined as
$$
\vcenter{\openup1\jot\halign{
\hfil$\displaystyle{#}$
&$\displaystyle{#}$\hfil
&$\displaystyle{#}$\hfil
\cr
\mu_{0,s}
&=1\quad
&\hbox{for $s\gt 0$}\ ,
\cr
\mu_{r,sr+m} 
&= \mu_{r-1,r}^s \mu_{m-1,m} +
    \sum_{p=0}^{s-1} \mu_{r-1,r}^p \mu_{r-1,(s-p)r+m}\hidewidth
\cr
&& \hbox{ for $r\ge 2,\,s\ge 1,\, 1\le m\lt r$}\ .
\cr
\mu_{r,sr} 
&=\sum_{p=0}^{s-1} \mu_{r-1,r}^p \mu_{r-1,(s-p)r}
   \qquad
& \hbox{ for $r\ge 1,\, s\ge 2$}\ .
\cr 
}}
\formula{museq.1}
$$
This looks quite different from~\frmref{museq}.  However, I claim
that~\frmref{museq} is just a harmless extension of~\frmref{museq.1}.
Indeed, just redefine the indexes by replacing $sr+m$ with $s$, also
accepting $s\ge 0$, which removes the implicit restriction $s\gt r$.
This is harmless, because for $s\le r$ one gets for~\frmref{museq}
$\mu_{r,s} = \mu_{r-1,s}$.  Therefore, in the second line one can
replace $\mu_{m-1,m}=\mu_{r-1,m}$ and include it into the sum.  This
completes the proof of lemma~\lemref{nrm.1}.\endproof

\subsection{3.4}{Completion of the proof of lemma~\lemref{generatrici}}
The statement of the lemma concerns only the sequence of generating
functions, that are estimated by~\frmref{nrm.1}.  The completion of
the proof rests on a suitable choice of the sequence $\{d_r\}_{r\ge
1}$, that was left arbitrary, and on a suitable estimate of the
sequence $\{\mu_{r,s}\}_{s\ge r\ge 0}$.  As a matter of fact, only the
diagonal elements of the latter sequence need to be estimated, because
the estimate for the generating functions in lemma~\lemref{nrm.1}
involves only $\mu_{r-1,r}=\mu_{r,r}$.

First, pick a value for $d$, with $0\lt d\lt 1/2$, and set
$$
d_r = \frac{b}{r^2}\ ,\quad b=\frac{6d}{\pi^2}\ ,
$$
so that $\sum_{r\ge 1}d_r=d$ in view of $\sum_{r\ge 1}1/r^2=\pi^2/6$.
The immediate consequence is that
$$
T_{r,s} \le \left(\frac{16}{b^2}\right)^{s-1}\ .
\formula{nrm.11}
$$
For, use the definition~\frmref{T} and recall the
definition~\frmref{jset} of $\Jscr_{r,s}$; then let $J\in\Jscr_{r,s}$
and evaluate
$$
\prod_{j\in J} \frac{1}{d_j} \le
 \frac{1}{b^{2(s-1)}} \prod_{j\in J} j^2\ ,
$$
because $\#(J)\le 2(s-1)$.  On the other hand one has
$$
\lgd \prod_{j\in J} j^2 =
  2 \sum_{j\in J} \lgd j \le 4(s-1)\ .
$$ 
This proves~\frmref{nrm.11}

Coming to the sequence~\frmref{museq}, the problem is to show that
$\mu_{r-1,r}\le \eta^{r-1}$ for some positive $\eta$.  For, only
$\mu_{r-1,r}$ enters the estimate~\frmref{nrm.1}.  By separating the
term $p=0$ in the sum one gets
$$
\eqalign{
\mu_{r,s}
&=\mu_{r-1,s} + 
   \mu_{r-1,r} \sum_{0\le q\lt s-r} \mu_{r-1,r}^q\mu_{r-1,s-r-rq}
\cr
&=\mu_{r-1,s} + 
   \mu_{r-1,r} \mu_{r,s-r}\ .
\cr
}
$$
Putting $r-1$ in place of $r$ and $r$ in place of $s$ in the latter
formula one gets
$$
\eqalign{
\mu_{r-1,r}
& = \mu_{r-2,r}+\mu_{r-2,r-1}\mu_{r-1,1}
\cr
& = \mu_{r-3,r}+\mu_{r-3,r-2}\mu_{r-2,2}+
       \mu_{r-2,r-1}\mu_{r-1,1}
\cr
&\ \ldots
\cr
& \le \mu_{0,r} + \mu_{0,1}\mu_{1,r-1} + \ldots
       \mu_{r-2,r-1}\mu_{r-1,1}
\cr
&\le \mu_{0,1}\mu_{r-2,r-1} + \ldots
      + \mu_{r-2,r-1}\mu_{0,1}\ .
\cr
}
$$
The last inequality requires a justification.  Just use
$$
\mu_{0,r}\lt\mu_{1,r}\lt\ldots\lt\mu_{r-1,r} =\mu_{r,r}=\ldots\ ,
$$ 
which is an immediate consequence of the definition.  Then, for
$r-j\ge j$ substitute $\mu_{r-j,j}=\mu_{j-1,j}$, and for $r-j\lt j$
use $\mu_{r-j,j}\lt\mu_{j-1,j}$, and so also
$\mu_{r-j,j}+1\le\mu_{j-1,j}$ which gets rid of the extra term
$\mu_{0,r}=1$.

Thus, the sequence $\{\nu_r\}_{r\ge 1}$ defined as 
$$
\nu_1=1\ ,\quad 
\nu_r
=\sum_{j=1}^{r-1} \nu_j\nu_{r-j}
 \quad\hbox{for $r\gt 1$}
\formula{nuseq}
$$
is a majorant of $\{\mu_{r-1,r}\}_{r\ge 1}$.  This is  known as the
Catalan's sequence, and one has
$$
\nu_r
=\frac{2^{r-1}(2r-3)!!}{r!}
\le 4^{r-1}\ ,
\formula{nrm.16}
$$
where the common notation $(2n+1)!!=1\cdot 3\cdot\ldots\cdot(2n+1)$
has been used.   

Thus, we conclude that $\mu_{r-1,r}\le 4^{r-1}$.  Inserting the latter
inequality and~\frmref{nrm.11} in~\frmref{nrm.1} the statement of
lemma~\lemref{generatrici} follows.

\section{4}{Proof of theorem~\thrref{nrm.1}}
Having established the estimate of lemma~\lemref{generatrici} on the
sequence of generating functions it is now a standard matter to
complete the proof of theorem~\thrref{nrm.1}.  Hence this section will
be less detailed with respect to the previous ones.

The situation to be dealt with is the following.  An infinite sequence
$\{\chi_r\}_{r\ge 1}$ of generating functions is given, with
$\chi_r\in\Pscr_{r+2}$ (a homogeneous polynomial of degree $r+2$)
satisfying $\norma{\chi_r}_{R}\le\beta^{r-1}G$ for some real vector
$R$ with positive components and some positive $\beta$ and $G$.
Define a corresponding sequence of canonical transformations
$(x^{(r-1)},y^{(r-1)}) =
\exp(\lie{\chi_r})(x^{(r)},y^{(r)})$.  By
composition one also constructs a sequence
$\{\Cscr^{(r)}\}_{r\ge 0}$ of canonical transformations
$(x^{(0)},y^{(0)}) = \Cscr^{(r)}(x^{(r)},y^{(r)})$
recursively defined as
$$
\Cscr^{(0)}= {\rm Id}\ ,\quad \Cscr^{(r)} = 
 \exp(\lie{\chi_r}) \circ
  \Cscr^{(r-1)}\ ,
$$
Id denoting the identity operator.  The problem is to prove the following
statements.

\item{(i)}Every near the identity canonical transformation defined via
the exponential operator $\exp(\lie{\chi_r})$ is expressed
as a power series which is convergent in a polydisk $\Delta_{\rho R}$
for some positive $\rho$.

\item{(ii)}For any function $f(x^{(r-1)},y^{(r-1)})$ analytic in
$\Delta_{\rho R}$ the transformed function is analytic in the same
polydisk, and moreover
$$
\Bigl. f(x^{(r-1)},y^{(r-1)})\Bigr|_{(x^{(r-1)},y^{(r-1)}) =
\exp(\lie{\chi_r})(x^{(r)},y^{(r)})}
=\bigl[\exp(\lie{\chi_r}) f
\bigr](x^{(r)},y^{(r)})\ .
$$

\item{(iii)}The sequence $\bigl\{\Cscr^{(r)}\bigr\}_{r\ge 0}$ of
canonical transformations converges for $r\to\infty$ to a canonical
transformation $\Cscr^{(\infty)}$ which is analytic in a polydisk
$\Delta_{(1-d)\rho R}$ for some positive $d\lt 1/2$.

\item{(iv)}For any function $f$ analytic in $\Delta_{\rho R}$ the
sequence recursively defined as $f^{(0)}=f$,
$f^{(r)}=\exp(\lie{\chi_r}) f^{(r-1)}$ converges for $r\to\infty$ to a
function $f^{(\infty)}$ that is analytic in $\Delta_{(1-d)\rho R}$, and
moreover one has
$$
f^{(\infty)} = f\circ \Cscr^{(\infty)}\ .
$$

\noindent
The statement~(i) actually reduces to Cauchy's proof of the existence
and uniqueness of the local solution of an analytic system of
differential equations.  For, the transformation defined via the
exponential operator is the time--one canonical flow induced by the
Hamiltonian vector field generated by $\chi_r$.  The statement~(ii)
actually claims that the substitution of variables in a function $f$
may be effectively replaced by the application of the exponential
operator to $f$; this is indeed the basis of the algorithm for
constructing the normal form used in sect.~\secref{2}.  A detailed
proof of both these statements may be found, e.g.,
in~\dbiref{Groebner-67}; however, the reader may be able to reconstruct
the proof by following the hints in~\dbiref{Giorgilli-95}.

The proof of~(iii) rests on the following remarks.  In the polydisk
$\Delta_{\rho R}$ one has $|\chi_r(x,y)| \le
\rho^{r+2}\norma{\chi_{r}}_{\rho R}$; this, in turn, implies that
$\bigl| x^{(r)}-x^{(r-1)} \bigr| \sim \beta^{r-1}\rho^{r+2}$ and
$\bigl| y^{(r)}-y^{(r-1)} \bigr| \sim \beta^{r-1}\rho^{r+2}$.  The
geometric bound on the latter quantities implies that $\sum_{r\gt
1} \bigl| x^{(r)}-x^{(r-1)} \bigr|$ and $\sum_{r\gt 1} \bigl|
y^{(r)}-y^{(r-1)} \bigr|$ behave like geometric series, i.e., converge
for $\rho$ small enough.  Thus, the claim follows from Weierstrass
theorem.  Finally, the statement~(iv) follows from~(ii) being true for
all $r\gt 0$, which implies that both sequences $f^{(r)}=\Cscr^{(r)}f$
and $f\circ\Cscr^{(r)}$ converge to the same limit.  This concludes
the proof of theorem~\thrref{nrm.1}.

\appendix{lyap.app.a}{Justification of the normalization algorithm}
Justifying the normalization algorithm of sect.~\secref{2} is just
matter of rearranging terms in the expansion of $\exp(\lie{\chi_r})
H^{(r-1)}$.  Considering first $H_0$ and $H^{(r-1)}_r$ together, one
has
$$
\vcenter{\openup1\jot\halign{
\hfil$\displaystyle{#}$
&$\displaystyle{#}$\hfil
&$\displaystyle{#}$\hfil
&$\displaystyle{#}$\hfil
\cr
\exp(\lie{\chi_r}) \bigl(H_0+H^{(r-1)}_r\bigr)
&= H_0\ &+ \lie{\chi_r}H_0\ &+\sum_{s\ge 2}
                               \frac{1}{s!}\lie{\chi_r}^s H_0
\cr
&&+H^{(r-1)}_r &+\sum_{s\ge 1}
                  \frac{1}{s!}\lie{\chi_r}^s H^{(r-1)}_r\ .
\cr
}}
$$
Here, $H_0$ is the first term in the transformed Hamiltonian $H^{(r)}$
in~\frmref{alg.1}.  In view of~\frmref{alg.2} one has
$\lie{\chi_r}H_0+H^{(r-1)}_r=Z_r$, which kills the unwanted term
$H^{(r-1)}_r$ and replaces it with the normalized term $Z_r$.  The two
sums may be collected and simplified by calculating
$$
\eqalign{
\sum_{s\ge 2}\frac{1}{s!}\lie{\chi_r}^s & H_0 +
 \sum_{s\ge 1} \frac{1}{s!}\lie{\chi_r}^s H^{(r-1)}_r
\cr
& =\sum_{s\ge 2}
    \frac{1}{(s-1)!}\lie{\chi_r}^{s-1}
     \left[\frac{1}{s}\left(\lie{\chi_r}H_0 + H^{(r-1)}_r
      \right) +
       \frac{s-1}{s} H^{(r-1)}_r
        \right]
\cr
& =\sum_{s\ge 2}
    \frac{1}{(s-1)!}\lie{\chi_r}^{s-1}
     \left(\frac{1}{s}Z_r +
      \frac{s-1}{s} H^{(r-1)}_r
       \right)\ .
\cr
}
$$
Here, both $\lie{\chi_r}^{s-1}Z_r$ and $\lie{\chi_r}^{s-1} H^{(r-1)}_{r}$
are homogeneous polynomials of degree $sr+2$, that are added to
$H^{(r)}_{sr}$ in the second of~\frmref{alg.3}.

Proceed now by transforming the functions $Z_1,\ldots,Z_{r-1}$
that are already in normal form.  Recall that no such term exists for
$r=1$.  For $r\gt 1$ calculate
$$
\exp(\lie{\chi_r}) Z_m = Z_m + 
 \sum_{s\ge 1} \frac{1}{s!}\lie{\chi_r}^s Z_m
  \ ,\quad
   \hbox{for $1\le m\lt r$}\ .
$$
The term $Z_m$ is copied into $H^{(r)}$ in~\frmref{alg.1}. The term
$\lie{\chi_r}^s Z_m$ is a homogeneous polynomial of degree $sr+m+2$
that is added to $H^{(r)}_{sr+m}$ in the first of~\frmref{alg.3}.

Finally, consider all terms $H^{(r-1)}_{s}$ with $s\gt r$, that may be
written as $H^{(r-1)}_{lr+m}$ with $l\ge 1$ and $0\le m\lt r$, the
case $l=1\,,\>m=0$ being excluded. One gets
$$
\exp(\lie{\chi_r}) H^{(r-1)}_{lr+m} = 
 \sum_{p\ge 0} \frac{1}{p!}\lie{\chi_r}^p H^{(r-1)}_{lr+m}
$$
where $\lie{\chi_r}^p H^{(r-1)}_{lr+m}$ is a homogeneous polynomial of
degree $(p+l)r+m$.  Collecting all homogeneous terms with $m=0$, $l\ge
2$ and $p+l=s\ge 2$ one gets $\sum_{p=0}^{s-2} \frac{1}{p!}\lie{\chi_r}^p
H^{(r-1)}_{(s-p)r}$, that is added to $H^{(r)}_{sr}$ in the second
of~\frmref{alg.3}.  Similarly, collecting all homogeneous terms with
$0\lt m\lt r$, $\>l\ge 1$ and $p+l=s\ge 1$ one gets $\sum_{p=0}^{s-1}
\frac{1}{p!}\lie{\chi_r}^p H^{(r-1)}_{(s-p)r+m}$, that is added to
$H^{(r)}_{sr+m}$ in the first of~\frmref{alg.3}.  The latter case
does not occur for $r=1$.  This completes the justification of the
formal algorithm.

\appendix{lyap.app.b}{Technical calculations}
The aim is to check the
estimates~\frmref{cauchy.f}, \frmref{cauchy.f-beq}
and~\frmref{cauchy.Z-dies}.   Write, generically, $\chi=\sum_{j,k}
c_{j,k}x^jy^k$ and $f=\sum_{j,k} f_{j,k}x^jy^k$. Then compute
$$
\lie{\chi}f = \sum_{j,k,j',k'} \sum_{l=1}^{n}
 \frac{j'_lk_l-j_lk'_l}{x_l y_l} c_{j,k}f_{j',k'} x^{j+j'}y^{k+k'}\ .
\formula{par_poisson}
$$
Using the definition of norm evaluate
$$
\eqalign{
\ordnorma{\lie{\chi}f}_{1-\delta} 
&\le 
  \sum_{j,k,j',k'} \sum_{l=1}^{n} 
   \frac{\left|j'_lk_l - j_lk'_l\right|}{R_l^2}
    |c_{j,k}|\, |f_{j',k'}| (1-\delta)^{|j+k|+|j'+k'|-2} R^{j+k+j'+k'}
\cr
&\le
  \frac{1}{\Lambda^2} \sum_{j,k} 
   \sum_{j',k'} 
    \sum_{l=1}^{n} \left|j'_lk_l - j_lk'_l\right|\,
     |c_{j,k}|\, \bigl((1-\delta')-(\delta-\delta')\bigr)^{|j+k|-1} R^{j+k}
\cr 
&\phantom{\Lambda^2 \sum_{l=1}^{n} {\left|j_lk'_l +
j'_lk_l\right|}\qquad}
 \times   
  |f_{j',k'}| \bigl((1-\delta'')-(\delta-\delta'')\bigr)^{|j'+k'|-1} R^{j'+k'}\ .
\cr
}
$$
If $f$ is a generic function, then 
in view of $|j_l|\le |j+k|$ and $k_l\le |j+k|$ one has 
$$
\sum_{l=1}^{n} |j'_lk_l-j_lk'_l| \lt |j+k|\sum_{l=1}^{n}|j'_l+k'_l| =
 |j+k|\cdot|j'+k'|\ .
\formula{ps.1}
$$
Replacing in the estimate above and using the elementary
inequality
$$
m(\lambda-x)^{m-1} \lt \frac{\lambda^m}{x}\quad {\rm for}\ 0\lt
x\lt\lambda\ {\rm and}\ m\ge 1
\formula{ineq}
$$
one gets
$$
\eqalign{
\ordnorma{\lie{\chi}f}_{1-\delta} 
&\le 
  \frac{1}{\Lambda^2}
   \sum_{j,k} 
    |c_{j,k}|\,|j+k|\bigl((1-\delta')-(\delta-\delta')\bigr)^{|j+k|-1} R^{j+k}
\cr
&\qquad\qquad
  \times
   \sum_{j',k'} 
    |f_{j',k'}|\, |j'+k'|\bigl((1-\delta'')-(\delta-\delta'')\bigr)^{|j'+k'|-1} R^{j'+k'}
\cr
&\le
  \frac{1}{(\delta-\delta')(\delta-\delta'')\Lambda^2}
   \sum_{j,k} |c_{j,k}|\,(1-\delta')^{|j+k|} R^{j+k}
\cr
& \phantom{\frac{1}{(\delta-\delta')(\delta-\delta'')\Lambda^2}}\qquad\qquad
   \times
    \sum_{j',k'} |f_{j',k'}| (1-\delta'')^{|j'+k'|} R^{j'+k'}\ ,
\cr
}
\formula{ps.3}
$$
from which~\frmref{cauchy.f} immediately follows in view of the
definition of the norm.  

In order to prove~\frmref{cauchy.f-beq} recall that $\chi^{\natural}$
contains only monomials $c_{j,k} x^j y^k$ with
$j+k\in\Kscr^{\natural}$.  The projection
$\bigl(\lie{\chi^{\natural}}f^{\flat}\bigr)^{\natural}$ is just part
of the general expression~\frmref{par_poisson}.  In particular the
value $l=1$ in the sum must be discarded because the resulting
monomials belong to $\Pscr^{\flat}$.  Moreover, for
$k\in\Kscr^{\natural}$ one has $\sum_{l=2}^n |j_l+k_l| \le 2$.  Thus
the estimate~\frmref{ps.1} may be replaced by
$$
\sum_{l=2}^{n} |j'_lk_l-j_lk'_l| \lt 2 \sum_{l=2}^{n}|j'_l+k'_l| =
 2 |j'+k'|\ .
\formula{ps.2}
$$
Hence the inequality~\frmref{ineq} must be used only for the term
involving $|j'+k'|$, and there is no need to introduce the divisor
$\delta-\delta'$. Use instead $1/(1-\delta) \lt 2$ in view of $\delta\lt 1/2$.

Coming finally to~\frmref{cauchy.Z-dies}, replace $f$ in the general
expression~\frmref{par_poisson} by
$Z^{\sharp}+Z^{\natural}=\sum_{\nu\in\Kscr^{\sharp}\cup\Kscr^{\sharp}} z_{\nu,\nu} x^{\nu}
y^{\nu}$.  Recall also that the coefficients $c_{j,k}$ of $\chi$ have
the form $c_{j,k}=\frac{\psi_{j,k}}{\langle k-j,\lambda\rangle}$, in
view of~\frmref{gen_soluzione}.  Then~\frmref{ps.1} may be replaced by
$$
\sum_{l} |\nu_l(j_l-k_l)| 
 \le |\nu|\sum_{l} |j_l-k_l|
  \le |\nu|\,|j-k|
\formula{ps.5}
$$
On the other hand, by lemma~\lemref{art.1}  one has 
$$
|c_{j,k}| \le \frac{|\psi_{j,k}|}{|j-k|\gamma}\ ,
$$
so that the factor $|j-k|$ in~\frmref{ps.5} is compensated by the
divisor here.  This removes the need to introduce the divisor $\delta-\delta'$
in the rest of the estimates. Use instead $(1-\delta)^{|j+k|-1} \le
2(1-\delta)^{|j+k|}$, which holds true in view of $\delta\lt 1/2$.
Then~\frmref{ps.3} is replaced by
$$
\displaylines{
\ordnorma{\lie{\chi}(Z^{\sharp}+Z^{\natural})}_{1-\delta} 
\le 
  \frac{1}{\Lambda^2}
   \sum_{j,k} 
    \frac{2}{\gamma}|\psi_{j,k}|\,(1-\delta)^{|j+k|} R^{j+k}
\hfill\cr\hfill
  \times
   \sum_{\nu\in\Kscr^{\sharp}} 
    |z_{\nu,\nu}|\, {|\nu|}\bigl((1-\delta'')-(\delta-\delta'')\bigr)^{2|\nu|-1} R^{2\nu}
\qquad\cr\hfill
\le
  \frac{1}{(\delta-\delta'')\gamma\Lambda^2}
   \sum_{j,k} |\psi_{j,k}|\,(1-\delta')^{|j+k|} R^{j+k}
    \sum_{\nu\in\Kscr^{\sharp}} |z_{\nu,\nu}| (1-\delta'')^{2|\nu|} R^{2\nu}\ .
\cr
}
$$
Thus~\frmref{cauchy.Z-dies} follows in view of the definition of the norm.

\references

\bye